\documentclass[reqno,11pt]{article}

\usepackage[mathscr]{eucal}
\usepackage[all]{xy}
\usepackage{mathrsfs}
\usepackage{xypic}
\usepackage{amsfonts}
\usepackage{amsmath}
\usepackage{amsthm}
\usepackage{amssymb}
\usepackage{latexsym}
\usepackage{eufrak}

\theoremstyle{plain}
\newtheorem{thm}[equation]{Theorem}
\newtheorem{prop}[equation]{Proposition}
\newtheorem{lem}[equation]{Lemma}
\newtheorem{cor}[equation]{Corollary}

\theoremstyle{definition}
\newtheorem{defn}[equation]{Definition}

\theoremstyle{remark}

\newtheorem{notation}[equation]{Notation}

\makeatletter
\renewcommand{\subsection}{\@startsection{subsection}{2}{0pt}{-3ex
plus -1ex minus -0.2ex}{-2mm plus -0pt minus
-2pt}{\normalfont\bfseries}} \makeatother

\numberwithin{equation}{subsection}

\newcommand{\erem}{\hphantom{.}\hfill$\lozenge$\end{rem}}

\newcommand{\bimod}[1]{#1\text{-}{\sf{bimod}}}

\newcommand{\hdot}{{\:\raisebox{3pt}{\text{\circle*{1.5}}}}}

\newcommand{\idot}{{\:\raisebox{2pt}{\text{\circle*{1.5}}}}}

\DeclareMathOperator{\Lie}{\mathrm{Lie}}
\DeclareMathOperator{\Tr}{\mathrm{Tr}}
\DeclareMathOperator{\Rep}{\mathrm{Rep}}

\DeclareMathOperator{\act}{\mathrm{act}}

\def\map{\longrightarrow}
\newcommand{\dis}{\displaystyle}

\newcommand{\beq}{\begin{equation}\label}
\newcommand{\eeq}{\end{equation}}

\newcommand{\iso}{{\;\stackrel{_\sim}{\to}\;}}
\newcommand{\cd}{\!\cdot\!}

\DeclareMathOperator{\GL}{\mathrm{GL}}

\def\ccirc{{{}_{\,{}^{^\circ}}}}

\newcommand{\g}[1]{\mathfrak{#1}}
\newcommand{\scr}[1]{\mathscr{#1}}

\DeclareMathOperator{\End}{\mathrm{End}}

\newcommand{\op}{\operatorname}
\DeclareMathOperator{\Hom}{\mathrm{Hom}}

\DeclareMathOperator{\Der}{\mathrm{Der}}
\DeclareMathOperator{\dder}{\mathbb{D}\mathbf{er}}

\newcommand{\exact}{_{\op{exact}}}
\newcommand{\closed}{_{\op{closed}}}

\newcommand{\De}{\Delta}

\newcommand{\id}{\mathrm{id}}
\newcommand{\Id}{\mathrm{Id}}

\newcommand{\im}{\mathrm{Im}}

\newcommand{\eex}{\hfill$\lozenge$\end{examp}}

\newcommand{\eps}{\varepsilon}

\newcommand{\cyclic}{_{\operatorname{cyclic}}}
\renewcommand{\a}{\!A}
\newcommand{\mcyclic}{_{\operatorname{anticyclic}}}
\newcommand{\bi}{\imath}
\newcommand{\bl}{{\scr L}}
\newcommand{\oma}{{\Om^1\a}}

\newcommand{\omb}{{\Om^\hdot\a}}

\newcommand{\rep}{{\Rep(A,V)}}

\newcommand{\bom}{{\overline{\Omega}}}
\newcommand{\HC}{{\overline{H\!C}}}
\renewcommand{\o}{{\otimes}}
\newcommand{\be}{\beta}

\DeclareMathOperator{\ad}{\mathrm{ad}}

\DeclareMathOperator{\Ker}{\mathrm{Ker}}
\DeclareMathOperator{\ev}{{\mathrm{ev}}}

\def\calT{{\mathscr{T}}}

\newcommand{\evom}{\ev_{_\Om}}
\newcommand{\BL}{{\mathbf{L}}}

\newcommand{\bbi}{{\mathbf{i}}}

\newcommand{\mto}{\longmapsto}

\newcommand{\inv}{^{-1}}
\newcommand{\vi}{${\sf {(i)}}\;$}
\newcommand{\vii}{${\sf {(ii)}}\;$}
\newcommand{\viii}{${\sf {(iii)}}\;$}

\newcommand{\sset}{\subset}

\newcommand{\into}{{}^{\,}\hookrightarrow^{\,}}
\newcommand{\too}{\,\longrightarrow\,}
\newcommand{\onto}{\twoheadrightarrow}
\newcommand{\tooo}{{\;{-\!\!\!-\!\!\!-\!\!\!-\!\!\!\longrightarrow}\;}}

\newcommand{\oper}{\operatorname}

\newcommand{\C}{\mathbb C}
\newcommand{\R}{\mathsf R}

\newcommand{\Z}{\mathbb Z}

\newcommand{\HH}{{H\!H}}

\newcommand{\DR}{{\oper{DR}}}
\newcommand{\starb}{\star_{_\beta}}

\renewcommand{\th}{\theta}

\newcommand{\munc}{{\mu_{\mathrm{nc}}}}

\newcommand{\Om}{\Omega}

\newcommand{\Th}{\Theta}

\def\ip<#1,#2>{\left\langle#1,#2\right\rangle}
\def\sp<#1>{\left\langle#1\right\rangle}

\newcommand{\BT}{{\mathsf{T}}}

\newcommand{\beps}{\C[\eps]/\langle\eps^2\rangle}

\def\ip<#1,#2>{\left\langle#1,#2\right\rangle}

\def\npb{\noindent$\bullet\quad$\parbox[t]{120mm}}

\newcommand{\om}{\omega}
\newcommand{\pa}{\partial}

\newcommand{\ka}{{\boldsymbol{\kappa}}}
\newcommand{\kt}{{\Bbbk\langle t\rangle}}
\newcommand{\al}{{\alpha}}

\newcommand{\en}{{\enspace}}
\newcommand{\br}{_\natural}
\newcommand{\wh}{\widehat}

\def\k{{\Bbbk}}

\begin{document}

\centerline{\huge{\textbf{\Large{Double derivations and Cyclic homology}}}}

\bigskip

\centerline{\large{\sc Victor Ginzburg}}
\medskip{\small
\begin{abstract}{\small We give a new construction of cyclic homology
of  an associative
algebra $A$ that does not involve Connes' differential. 
Our approach is based on an extended version of
the complex  $\omb,$ of  noncommutative
differential forms on $A,$ and  is similar in spirit to the
de Rham approach to equivariant cohomology. Indeed, our extended complex
maps naturally to the equivariant de Rham complex of any
representation scheme~$\Rep A$.

We define
cyclic homology as the
cohomology of the total complex $(\omb)[u], d+u\cdot\bi_\De$,
arising from two anti-commuting differentials, $d$ and $\bi_\De$,
on $\omb$. The  differential $d$,
that replaces the Connes differential $B$, is the Karoubi-de Rham
differential. The  differential $\bi_\Delta$,
that replaces the Hochschild differential $b$,
is a  map  analogous to contraction with
a vector field. This new map  has no
commutative counterpart.}
\end{abstract}
}

\section{Introduction}
Throughout, we fix a  field $\k$ of characteristic 0
and write $\o=\o_\k.$
 By an  algebra we will always mean
an  associative unital $\k$-algebra.

\subsection{}
It is well-known that a regular vector field on a smooth affine
algebraic variety $X$ is the same thing as a derivation
$\k[X]\to \k[X]$, of the coordinate ring of $X$. Thus,
derivations of a commutative algebra $A$ play the role
of vector fields.

It has been commonly accepted until recently that this
 point of view 
 applies to noncommutative algebras $A$ as well.
A first indication towards a different point of view
was a discovery by Crawley-Boevey \cite{CB} that,
for a smooth affine
curve $X$ with coordinate ring $A=\k[X]$,
the algebra of differential operators on $X$
 can be constructed by means of {\em double derivations}
$A\to A\o A$, rather than ordinary derivations
$A\to A.$ Since then, the significance
of double derivations in noncommutative geometry
was  explored further in \cite{VdB} and \cite{CBEG}.

To explain the role of  double derivations in more detail
we need to introduce, in \S\S\ref{dder}-\ref{free_prod},
 some basic definitions involving 
 various kinds of derivations and also 
 {\em free products}.

\subsection{Derivations and $t$-derivations.}\label{dder}
 Let $B$ be any algebra  and $M$ a $B$-{\em bimodule}.

Recall that a  $\k$-linear map $ F: B\to M$
is said to be a derivation of $B$ with coefficients in $M$
if  $ F(b_1b_2)= F(b_1)\,b_2+b_1 \,F(b_2),\,
\forall b_1,b_2\in B.$
Given  a subalgebra $R\sset B$,
 we let $\Der_R(B,M)$ denote the space of relative derivations 
of $B$ with respect to the subalgebra $R$, that is,
of derivations $B\to M$ that annihilate the subalgebra $R$.
 In the `absolute' case $R=\k$
we will use simplified notation $\Der(B,M).$

Now, fix  an element $t\in B$.
We will be often interested in 
 $\k$-linear maps $F: B\to M$ such that
\beq{tder}
F(b_1\, t\, b_2)=F(b_1)\, t\, b_2 + b_1\, t\,  F(b_2),
\quad\forall b_1,b_2\in B.
\eeq

\begin{defn}\label{tder_def} A  $\k$-linear map
$F: B\to M$ is called a $t$-{\em derivation} if equation
\eqref{tder} holds for $F$ and, moreover, we have $F(1)=0.$
\end{defn}

From equation \eqref{tder}, by induction on $n\geq 1$
one proves
$$F(t^n)=
t^n\, F(1)+ t^{n-1}\, F(1)\, t+ t^{n-2}\, F(1)\, t^2+\ldots+
F(1) \, t^n.
$$
It follows that, for any $t$-derivation $F$ and an
integer $n\geq 1$, one has $F(t^n)=0;$
also one has $F(t\, b)=t\, F(b)$ and
$F(b\, t)=F(b)\, t,$ for all $b\in B$.

Let $\kt$ be the subalgebra in $B$ generated by the element
$t$. It is clear that: 
{\em Any  derivation
$F\in \Der_\kt(B,M)$ is automatically a $t$-derivation}.

 We warn the reader that  although 
any $t$-derivation annihilates the subalgebra $\kt$,
a $t$-derivation  need not necessarily be a derivation,
i.e., an element of $\Der_\kt(B,M),$ in general.

\begin{lem}\label{tder_lem2} For $t\in B$, the space of $t$-derivations
$B\to B$ is a Lie algebra with respect to the commutator bracket.
\end{lem}
\proof For any $t$-derivations $F,G: B\to B$ and $b_1,b_2\in B,$
we compute
\begin{align*}
&GF(b_1\, t\, b_2)=
G(F(b_1)\, t\,  b_2 + b_1\,  t\, F(b_2))\\
&=
GF(b_1)\, t\,  b_2 +F(b_1)\, t\,  G(b_2)
+ G(b_1)\,  t\, F(b_2)
+b_1\,  t\, GF(b_2).
\end{align*}
Performing a similar computation for $FG,$ we find that  the commutator 
$(G F-F G)(b_1\, t\, b_2)=[G,F](b_1\, t\, b_2)$ equals
\begin{align*}
[G,F](b_1)\, t\,  b_2&+F(b_1)\, t\,  G(b_2)
+G(b_1)\,  t\, F(b_2)-G(b_1)\, t\,  F(b_2)
-F(b_1)\,  t\, G(b_2)\\
&+b_1\,  t\, [G,F](b_2)
=[G,F](b_1)\, t\,  b_2 +b_1\,  t\, [G,F](b_2).\qquad\quad\Box
\end{align*}

\subsection{Free product construction.}\label{free_prod}
Recall that a free product of two  algebras ${A}$ and $B$, 
is an associative algebra ${A}*B$ 
 whose elements are formal $\k$-linear combinations of words
$a_1b_1a_2b_2\ldots a_nb_n,$ for any $n\geq 1$ and
$a_1,\ldots,a_n\in {A},\,b_1,\ldots,b_n\in B.$
These words are taken up to equivalence  imposed
by the relation $1_{A}=1_{B}$; for instance, 
we have $\ldots b1_{A} b'\ldots=\ldots b1_{B}b'\ldots=\ldots
(b\cdot b')\ldots,$
for any $b,b'\in B.$ 

In the
special case where $B=\k[t]$, a polynomial algebra in one variable, we use
simplified notation ${A}_t:={A}*\k[t].$
It is convenient to view $A_t$ as
 an {\em even} graded algebra
 $A^\hdot_t=\oplus_{k\geq 0} A^{2k}_t$ with the grading
 that  counts {\em twice} the number of
occurences of the variable~$t$,
i.e., such that  $\deg t=2$
and  such that  the subalgebra $A\sset A_t$ is assigned grade degree
zero.
For each integer $k$, let $\Der_{\k[t]}^{2k}(A_t,A_t)$ 
denote the space  of graded
derivations $A^\hdot_t\to A^{\idot+2k}_t,$
 relative to the
subalgebra $\k[t]\sset A_t.$
Thus, we have a graded Lie algebra
$$\Der^\hdot_t(A_t)=\oplus_{k\in\Z}\Der_t^{2k}(A_t),\quad\text{where}
\quad \Der_t^{2k}(A_t):=\Der_{\k[t]}^{2k}(A_t,A_t).$$

View $A$ and $\k[t]$  as two subalgebras in $A_t$.

\begin{lem}\label{f_tlem}
Let $M$ be an $A_t$-bimodule, also viewed as an $A$-bimodule.

\noindent
\vi Restriction to the subalgebra $A\sset A_t$ provides a
vector space isomorphism
$$
\xymatrix{
{\left\{
{t\text{\em -derivations }}\atop
{F: A_t\to M}
\right\}}
\ar[rrr]_<>(0.5){\sim}^<>(0.5){F\mapsto f:= F|_A}&&&
{\left\{
{\k\text{\em -linear maps } f: A\to M}\atop
{\text{\em such that } f(1)=0}
\right\}}
}
$$
\vii The   isomorphism in \vi restricts to a  bijection:\hfill\break
\centerline{$\Der_{\k[t]}(A_t,M)\iso\Der(A, M).$}
\end{lem}

\begin{proof} It is clear that the assignment
$F\mapsto f:= F|_A$ gives an injective map
from the set of $t$-derivations $F: A_t\to M$
to the set of $\k$-linear maps $f: A\to M$
such that  $f(1)=0.$
We construct a map in the opposite direction
by assigning to any  $\k$-linear map $f: A\to M,$
such that  $f(1)=0,$ a map $f_t: A_t\to M$ given,
for any $a_1,\ldots,a_n\in A,$
by the following  
Leibniz type formala
\beq{f_t}
a_1\,t\,a_2\,t\ldots t\,a_n\mto
\sum\nolimits_{k=1}^n a_1\,t\ldots a_{k-1}\,t\,f(a_k)\,t\,
a_{k+1}\,t\ldots t\,a_n.
\eeq

One verifies that the map $f_t$
thus defined satisfies \eqref{tder}.
It is also immediate to check that the  maps
$F\mapsto F|_A$ and $f\mapsto f_t$
are inverse to each other.
 This proves (i).
Part (ii) is straightforward and is  left to the reader.
\end{proof}

\begin{notation}\label{f_tnot}
We write $f_t$ for the $t$-derivation
$A_t\to M$
given by formula \eqref{f_t}. Thus, $f_t$ is the $t$-derivation
corresponding to a $\k$-linear map $f: A\to M$ under
the bijection inverse to the isomorphism of Lemma~\ref{f_tlem}(i).
\hfill$\lozenge$
\end{notation}

For any derivation $\th: A\to A$
(viewed as a derivation $A\to  A_t$), the corresponding
 derivation $\th_t: A_t\to A_t$ clearly
has degree zero, i.e., we have
$\th_t\in \Der_t^0(A_t).$
\smallskip

We will  use symbolic 
Sweedler notation to write  a $\k$-linear map $\Th: A\to A\o A$,
as $a\mto \Th'(a)\otimes\Th''(a),$
where we systematically
suppress the summation symbol.
We always view
the space $A\otimes A$ as an $A$-bimodule with respect
to the {\em outer} bimodule structure, which is
defined by  the formula
$b(a'\otimes a'')c:=(ba')\otimes (a''c)$,
for any $a',a'',b,c\in A.$ 
\begin{notation} We  write $\dder\a:=\Der(A,A\o A)$ for the  corresponding
vector space of 
derivations, to be referred to as
{\em double derivations}.
\hfill$\lozenge$
\end{notation}

Next, to any  $\k$-linear map $\Th: A\to A\o A$, we associate 
 a  $\k$-linear map $A\to A_t$ given by
 the assignment $a\mto \Th'(a)\,t\,\Th''(a).$
Assuming further that $\Th(1)=0$, 
we may extend the latter map 
to get a  $t$-derivation  $\Th_t: A_t\to A_t$,
see Notation \ref{f_tnot} and Lemma \ref{f_tlem}.

It is clear that, for any 
double derivation $\Th\in\dder A$, the corresponding
derivation $\Th_t: A_t\to A_t$
has degree 2, i.e., we have
$\Th_t\in \Der_t^2(A_t).$

It is easy to check that  converse
statements are also true.
Specifically,  
using   a graded version of Lemma \ref{f_tlem}, 
one proves the following result that clarifies
the relationship between derivations
and $t$-derivations.

\begin{lem}\label{gradings} \vi The composite map 
\[
\xymatrix{
\Der(A,A)\ar[rr]^<>(0.5){A\into
A_t}&&\Der(A,A_t)\ar[rr]^<>(0.5){\th\mapsto \th_t}&&
\Der^\hdot_{\k[t]}(A_t,A_t)
}
\]
yields a Lie algebra isomorphism
$\Der(A,A)\iso\Der_t^0(A_t).$

\vii Similarly, the assignment $\Th\mapsto\Th_t$ yields a vector space 
 isomorphism
$\dder A\iso\Der_t^2(A_t).$\qed
\end{lem}

\subsection{The role of double derivations.}
To explain the role of double derivations, recall first that
ordinary  derivations of $A$ may also
be thought of,
algebraically,
as `infinitesimal automorphisms'. 
Specifically,  let $A[t]=A\o\k[t]$ be the
polynomial ring in one variable with coefficients in $A$.
This is a $\k[t]$-algebra and, for
any $\k$-linear map $\xi: A\to A$,
the assignment  $A\to A[t],$
$a\mapsto t\cdot \xi(a)$
can be 
uniquely extended  to  a 
$\k[t]$-linear map
 $t\xi: A[t]\to A[t]$. 

A well known elementary calculation yields

\begin{lem}\label{[t]}  The following properties of a $\k$-linear map  $\xi: A\to A$
are equivalent:

\npb{The map $\xi$ is a derivation of the algebra $A$;}

\npb{The map
$t\xi: A[t]\to A[t]$ is a  derivation of the algebra $A[t]$;}

\npb{The map
$\Id + t\xi:\
 A[t]/t^2\cd A[t]\to A[t]/t^2\cd A[t]
$
 is an
algebra automorphism.}
\end{lem}

All the above holds true, of course, no matter whether
the algebra $A$ is commutative or not.
Yet, the element $t$, playing
the role of formal parameter,
 is by definition a {\em central} element
of the algebra  $A[t]$.  

In noncommutative geometry, 
the assumption that the formal parameter be central
is not quite natural, however.
Thus, we are led to consider the  algebra $A_t=A*\k[t],$
freely generated by $A$ and an indeterminate $t$.

We are going to argue that, once the polynomial algebra
$A[t]$
is replaced by  $A_t$, the  free product algebra,
it becomes more natural to replace
derivations $A\to A$ by   double derivations
$A\to A\o A$. 

To see this, let $A_t^+=A_t\cd t\cd A_t$ be the two-sided ideal of the algebra
$A_t$ generated by  $t$. Thus, we have 
$A_t/A_t^+=A$. Further,  a free product analog of an obvious
isomorphism $A[t]/t^2\cd A[t]= A\oplus t\,A$
is  a natural $A$-bimodule isomorphism
\beq{AtA}
 A_t/(A^+_t)^2\stackrel{_\sim}\longleftrightarrow
A\oplus (A\o A),
\quad  a+ a'\,t\,a''\longleftrightarrow a \oplus (a'\o a'').
\eeq 

Thus, we see  that the  direct summand
$t\,A\sset  A[t]/t^2\cd A[t]$ gets replaced, in the free product setting,
by $A\o A$. Furthermore, with the notation of Lemma
\ref{gradings},
a free product analog of Lemma~\ref{[t]} reads.

\begin{lem}\label{del_der}  The following properties of a $\k$-linear map  
 $\Th: A\to A\o A$
are equivalent:

\npb{The map $\Th$ is a double derivation, i.e., $\Th\in\dder\a$;}


\npb{We have
$\Th_t\in \Der_t(A_t)$;}

\npb{The map
$\Id + \Th_t:\
 A_t/(A_t^+)^2\to A_t/(A_t^+)^2
$
 is an
algebra automorphism.}
\end{lem}

Thus, we conclude that in noncommutative geometry
the algebra  $A_t$
should play the role of  the polynomial algebra 
$A[t]$. Some aspects of this philosophy will be
discussed further in subsequent sections.

\subsection{Layout of the paper.}  In \S2, we recall the
definition of the DG algebra of noncommutative
differential forms. We also introduce the Karoubi-de Rham complex
and an extended  Karoubi-de Rham complex, that will
play a crucial role later.
In \S3 we develop the basics of Noncommutative Calculus
involving the action of double derivations on the
 extended  Karoubi-de Rham complex, via Lie derivative and
contraction operations.

In section 4, we construct a canonical morphism
from the extended de Rham complex of an algebra $A$
to the {\em equivariant} de Rham complex of
the corrresponding representation scheme $\Rep(A,V)$.
In \S5 we state two main theorems of the paper.
These theorems provide a description, in terms of the de Rham complex,
 of Hochschild homology
of an algebra $A$ and of cyclic homology of $A$, respectively.
The proofs of the  theorems are given in \S6.
The proofs are based on the properties of  the Karoubi operator
and on  Harmonic decomposition
of the algebra of  noncommutative
differential forms introduced by Cuntz-Quillen,
\cite{CQ1}, \cite{CQ2}.

The final section contains an alternative construction
of the extended de Rham complex. We also  outline
a new version of deformation theory of associative algebras
with {\em noncentral} deformation parameter.

\subsection{Acknowledgements.}{\footnotesize{  I am indebted to Boris Tsygan
for explaining to me his unpublished results \cite{T}.
I would like  to thank   
Vladimir Drinfeld and Michel Van den Bergh
for  useful discussions, and  Jim Stasheff for helpful comments.
This work  was partially supported by the NSF grant
DMS-0303465 and CRDF grant {{RM1-2545-MO-03.}}}}

\section{The extended Karoubi-de Rham complex}
\subsection{The commutator quotient.} 
Let $B=\oplus_{k\in \Z} B^k$ be a $\Z$-graded
algebra and $M=\oplus_{k\in \Z} M^k$ a graded $B$-bimodule.
A linear map $f: B^\hdot\to M^{\idot+n}$ is said
to be a {\em degree} $n$  {\em graded} derivation if, for any
$u\in B^k$ and $v\in B,$ we have
$f(uv)=f(u)\cdot v+ (-1)^{nk}u\cdot f(v)$. We write $\Der^n(B,M)$ for the vector
space of degree $n$ graded derivations. The direct sum 
$\Der^\hdot B:=
\bigoplus_{n\in \Z}\Der^n(B,B),$ of graded derivations of the
algebra $B$, has a natural Lie {\em super}-algebra structure
given by the {\em super}-commutator.

We write
$[B,B]$ for the  {\em super}-commutator space of a graded algebra $B$,
the $\k$-linear span of the
set  $\{uv-(-1)^{pq} vu\mid$
$u\in B^p,v\in B^q,\,p,q\in\Z\}.$ This
is a graded subspace of $B$, and we may consider
the commutator quotient  space $\R B:=B/[B,B]$,
equipped with induced grading $\R^\hdot B=\oplus_{k\in \Z} \R^kB.$
Any degree $n$ graded derivation  $f: B^\hdot\to B^{\idot+n}$ descends to a well-defined
linear map of graded vector spaces $\R(f): \R^\hdot B\to\R^{\idot+n}B$.

\begin{notation} 
Given an algebra $A$ and an
 $A$-bimodule  $M$, let 
 $T_A M=\oplus_{n\geq 0}T^n_AM$ be the tensor algebra of 
$M$ over $A$. Thus, $T_A^\hdot M$ is a graded associative
algebra with $T^0_AM=A$. 
In the special case $A=\k$ we will use unadorned notation
$T M:=T_\k M,$
and write $M^{\o n}:=T^nM.$
\end{notation}

Let  $\BT: M^{\o n}\to M^{\o n},\,
u_1\o u_2\o\ldots \o
u_n\to u_n\o u_1\o u_2\o\ldots \o u_{n-1}$
denote the generator of the group $\Z/n\Z$ of {\em cyclic permutations}
 of the tensor factors, and write
 $M^{\o n}\cyclic:=M^{\o n}/(\Z/n\Z)$ for the space of $\Z/n\Z$-{\em
coinvariants}.
It is clear that, for any $n=1,2,\ldots,$ we have  natural identifications
\beq{cyclic_word}
\R^n(T M)= M^{\o n}\cyclic=M^{\o n}/\im(\Id-\BT),
\eeq
for instance,
modulo commutators, one has
$x\o y\o z=z\o x\o y=y\o z\o x.$

Later on, we will also use the following space
\beq{sign}
M^{\o n}\mcyclic:=M^{\o n}/\im(\Id-(-1)^{n-1} \BT),\quad\forall n=1,2,\ldots,
\eeq
 of
coinvariants with respect to a {\em sign-twisted}
$\Z/n\Z$-action on $M^{\o n}$. Note that $\BT$ is a permutation
of parity $(-1)^{n-1}$. 
\subsection{Commutator quotient for free products.}\label{ttau}
We will be using free product constructions for {\em graded}
algebras. For a graded algebra $B$, the  algebra $B_t=B*\k[t]$
acquires  a natural grading $B^\hdot_t=\oplus_{k\in\Z} B^k_t,$
that agrees with the one
on the subalgebra $B\sset B_t$ and such that $\deg t=2$. 
Let $\R^\hdot(B_t)=\oplus_{k\in\Z} \R^k(B_t)$
be  the corresponding commutator quotient.

There is  a graded analog of Lemma \ref{f_tlem}.
Specifically, given a graded $B^\hdot_t$-bimodule
$M^\hdot=\oplus_{k\in\Z} M^k$ and
a  $\k$-linear map $f: B^\hdot\to M^{\idot+n}$ such that $f(1)=0$,
one introduces a degree $n$ derivation $f_t: B_t^\hdot\to M^{\idot+n}$
defined, for any homogeneous elements $u_1,\ldots,u_n\in B,$ by 
the formula
\begin{align}\label{f^-}
f_t&(u_1\,t\,u_2\,t\ldots t\,u_r)\\
&=
\sum\nolimits_{k=1}^r (-1)^{n(\deg u_1+\ldots+\deg u_{k-1})}\cd
u_1\,t\ldots u_{k-1}\,t\,f(u_k)\,t\,
u_{k+1}\,t\ldots t\,u_r.\nonumber
\end{align}

We will also
consider an {\em odd} version of free product.
 Specifically, let
 $B_\tau:=B*\k[\tau]$, and  equip this free product with a
$\Z$-grading  $B_\tau=\oplus_{k\in\Z} B^k_\tau$  that agrees with the
one on
$B$  and such that $\deg\tau=-1$ (as opposed to $\deg
t=2$). Thus, $B_\tau$ is isomorphic to $B_t$ as an associative
algebra but not as a graded algebra.

The grading on $B_t$, resp. on $B_\tau$,
induces a
grading $\R^\hdot(B_t)=\oplus_{k\in \Z} \R^{k}(B_t),$
 resp.  $\R^\hdot(B_\tau)=\oplus_{k\in\Z} \R^k(B_\tau).$

Now, given an algebra $B$ {\em without} grading, we may
(and will) view it as a graded algebra concentrated in degree zero.
In that case, the
resulting grading $B_t=\oplus_{k\geq 0} B^{2k}_t$ is the
one considered in \S\ref{free_prod}, that is, the even grading that
counts twice the number of occurencies of the variable $t$.
Similarly, we have the algebra $B_\tau=\oplus_{k\leq 0} B^k_\tau,$
with {\em non-positive} grading, and also,  graded
 commutator quotient spaces
 $\R^\hdot(B_t)=\oplus_{k\geq 0} \R^{2k}(B_t),$
 resp.  $\R^\hdot(B_\tau)=\oplus_{k\leq 0} \R^k(B_\tau).$

In particular, we have $\R^0(B_t)= \R^0(B_\tau)=\R B=B/[B,B],$
where
$[B,B]\sset B$ is the subspace spanned by ordinary commutators.
For any $n\geq 1,$ the space $\R^{2n}(B_t)$, resp.
the space $\R^{-n}(B_\tau),$ is spanned
by {\em cyclic} words $u_1\,t\,u_2\,t\ldots t\,u_n\,t,$
for instance, modulo (graded) commutators, we have
$u_1\,t\,u_2\,t=u_2\,t\,u_1\,t,$ resp.
$u_1\,\tau\,u_2\,\tau=
-u_2\,\tau\,u_1\,\tau.
$

This way, one arrives at
the following result that generalizes
\eqref{AtA}. 

\begin{lem}\label{triv} For  any  (non-graded) algebra $B$
and $n\geq 1$,
the assignment
$u_1\,t\,u_2\,t\ldots t\,u_n\mto$
$u_1\o u_2\o\ldots \o u_n$ yields 
natural   vector
space  isomorphisms
$B^{2n}_t\cong B^{\o n}\cong B^{-n}_\tau.$

 The induced isomorphisms of the corresponding
graded commutator quotients read
$$
\R^{2n}(B_t)\cong B^{\o (n-1)}\cyclic,\quad\text{\em resp.},\quad
\R^{-n}(B_\tau)\cong B^{\o (n-1)}\mcyclic,\quad\forall n\geq 1.
$$

For $n=0,1,$ we have 
$$\R^0(B_t)=B/[B,B]=\R^0(B_\tau),\quad\text{\em and}\quad
\R^2(B_t)\cong B\cong
\R^{-1}(B_\tau).
$$
\end{lem}

Observe next that
the assignment
$\frac{d}{d\tau}:\tau\mapsto 1,\,b\mapsto 0,\,\forall b\in B,$
extends uniquely to a degree $1$ graded derivation
$\frac{d}{d\tau}: B_\tau\to B_\tau$. In particular, we have
$$\mbox{$\frac{d}{d\tau}$}(\tau^n)=1\cdot\tau\cdot\ldots\cdot\tau-
\tau\cdot 1\cdot\tau\cdot\ldots\cdot\tau+\ldots+ (-1)^{n-1}
\tau\cdot\ldots\cdot\tau\cdot 1=
\begin{cases} 0&\text{if}\; n\;\text{is even}\\
\tau^{n-1}&\text{if}\; n\;\text{is odd.}
\end{cases}
$$
It follows easily that one has
  $\frac{d}{d\tau}\ccirc
\frac{d}{d\tau}=0$,
i.e.,
the map $\frac{d}{d\tau}:B^\idot_\tau\to B^{\idot+1}_\tau$ is a differential.


Recall further that, associated with any algebra $B$,
one has a standard {\em bar complex}
\beq{bar}
\ldots\stackrel{\pa}\too B^{\o 4}
\stackrel{\pa}\too B^{\o 3} \stackrel{\pa}\too B^{\o 2}
\stackrel{\pa}\too B.
\eeq
Furthermore, it is well known that the bar differential descends
to  well defined maps, cf. \cite{Lo}
\beq{cyc_comp}
\ldots\stackrel{\pa}\too B^{\o 4}\mcyclic\stackrel{\pa}\too B^{\o 3}\mcyclic
\stackrel{\pa}\too B^{\o 2}\mcyclic\stackrel{\pa}\too B.
\eeq
The resulting complex is known as the {\em cyclic complex} for $B$.

Verification of the following result is straightforward.

\begin{prop}\label{drin}
\vi The  bar differential $\pa$ in \eqref{bar} is transported, via
the   isomorphisms $B^{\o n}\cong B^{-n}_\tau$
of Lemma \ref{triv}(i),
to the  differential $\frac{d}{d\tau}:{B^\idot_\tau\to~B^{\idot+1}_\tau.}$

\vii The induced isomorphism
of Lemma \ref{triv}(ii), provides
a natural identification of the cyclic complex \eqref{cyc_comp}
with the 
complex
$\bigl(\R^{<0}(B_\tau),\,\R(\frac{d}{d\tau})\bigr).$
\end{prop}

I first learned part (i) of the above result, which is a special
case $M=0$ of Proposition \ref{tsy} of \S7 below, from V. Drinfeld.

\subsection{Noncommutative differential forms.}
We  fix  an 
algebra $B$, a subalgebra $R\sset B$, and
 write
$B\otimes_R B\stackrel{m}\to B$ for the multiplication map.
Let $\Om^1_RB:=\Ker(m)$
 be  the $B$-bimodule of noncommutative
  1-forms on $B$ relative to the subalgebra $R$, see
  \cite[Sect.~2]{CQ1},
and let  
 $i_\Delta: \Om^1_RB\into B\otimes_R B$ denote the tautological 
imbedding (this notation will be justified later).
We have a
fundamental short exact sequence of $B$-bimodules
\beq{fund}
0\map\Om^1_RB\stackrel{i_\De}\too B\o_R B\stackrel{m}\too B\map 0.
\eeq

The assignment $b\mapsto db:=1\o b-b\o 1$ gives
 a canonical
derivation $d:B\to \Om^1_RB.$ This derivation is `universal'
in the sense
that, for any $B$-bimodule $M$, we have a bijection
\beq{der_def}
\Der_R(B,M)\iso\Hom_{\bimod{B}}(\Om^1_RB, M),\quad
\th\mapsto i_\th,
\eeq
where the map $i_\th:\Om^1_RB\to M$ is defined by
the formula $i_\th(u\,dv):=u\cdot\th(v).$

The tensor algebra
  $\Om_R^\hdot B:=T_B^\hdot(\Om^1_RB),$ of the 
$B$-bimodule  $\Om^1_RB,$ is a DG algebra  $(\Om^\hdot_RB ,d),$
called
the  algebra of noncommutative  differential forms on $B$
relative to the subalgebra $R$
(we will interchangeably use the notation $\Om_R B$ or
$\Om^\hdot_RB $ depending on whether we want to emphasize the grading or
not).
For each $n\geq 1,$ there is a standard
isomorphism of left $B$-modules,
see \cite{CQ1}, $\Om^n_RB=B\o_R T_R^n(B/R);$ usually, one writes
$b_0\,db_1\,db_2\ldots db_n\in \Om^n_RB$ for the $n$-form
corresponding to an element $b_0\o(b_1\o\ldots\o b_n)\in
B\otimes_R T_R^n(B/R)$ under this isomorphism.
The de Rham differential  $d: \Om^\hdot_RB \to\Om^{\idot+1}_RB$
is defined by the formula
$d: b_0\,db_1\,db_2\ldots db_n\mapsto
db_0\,db_1\,db_2\ldots db_n.$

Following Karoubi, we
define the noncommutative de Rham complex of
$B$ (relative to the subalgebra $R$)
to be  the commutator quotient space 
$$\DR_R  B:=\R(\Om_R B)=\Om_R B/[\Om_R B,\Om_R B],$$
 of the
{\em graded} algebra $\Om^\hdot_RB $. The space $\DR_R B$ comes
equipped with a natural grading and
with the de Rham differential $d: \DR^\hdot_R B\to \DR_R^{\idot+1}B$,
induced from the one on $\Om^\hdot_RB $,
c.f. also \cite{Lo} for more details.
In degree zero, we have $\DR^0_RB=\R B=B/[B,B]$.

In the `absolute' case $R=\k$ we will use unadorned notation
$\Om^nB:=\Om^n_\k B,
\,\DR B:=\DR_\k B,$ etc.

Let $\beps$ be the ring of dual numbers. Equip the free
product
$B_\eps:=B*(\beps)$  with a grading such that
$\deg B=0$ and $\deg\eps=1$. The graded commutator map
$\ad\eps: u\mapsto [\eps,u]:=\eps\,u-u\,\eps$
gives a degree 1 derivation $\ad\eps: B_\eps^\idot\to
B_\eps^{\idot+1},$ moreover, we have
$(\ad\eps)^2=\ad(\eps^2)=0$. Thus,
$(B_\eps,\,\ad\eps)$ is a DG algebra.

Further, the  derivation $B\to B_\eps,\,b\mapsto [\eps,b]$ gives rise,
via the universal property of $\Om^\hdot B$, to
a  DG algebra imbedding
$(\Om^\hdot B,\,d)\into (B_\eps,\,\ad\eps)$, cf. \cite[Proposition
1.4]{CQ1},
 given by the formula
\begin{align*}
b_0\,db_1\,db_2\ldots db_n&\mto
b_0\,[\eps,b_1]\,[\eps,b_2]\ldots [\eps,b_n]\\
&\;=\;
b_0\,\eps\, b_1\,\eps\, b_2\ldots \eps\, b_n\,+\,
(-1)^n\, b_0\,b_1\,\eps\, b_2\ldots \eps\, b_n\,\eps.
\end{align*}

\subsection{The extended de Rham complex.} Fix an algebra $A$.
We introduce an enlargement of the
noncommutative de Rham complex,  $\DR A,$ as follows. 

First, form the free product algebra $A_t=A*\k[t]$.
This free product contains $\k[t]$ as a subalgebra
and we consider 
$\Om_t A:=\Om_{\k[t]}(A_t)$,
the DG algebra of noncommutative differential forms on $A_t$
relative to the subalgebra $\k[t]$.
We define  extended
de Rham complex of $A$ as
$$\DR_t A:=\DR_{\k[t]}(A_t)=\Om_{\k[t]}(A_t)/[\Om_{\k[t]}(A_t),
\Om_{\k[t]}(A_t)].
$$

On the other hand, one may form a free product $(\Om A)_t=(\Om A)*\k[t],$
and equip it with the standard even grading that counts twice 
the number of occurencies of $t$ (and disregards the degrees of differential
forms). We have the corresponding
 direct sum decomposition into graded components
\beq{dsd}(\Om A)_t/[(\Om A)_t,(\Om A)_t]=\R((\Om A)_t)=\oplus_{k\geq 0}
\R^{2k}((\Om A)_t).
\eeq
The differential
$d: \Om^{\,}\a  \to \Om^{\,}\a\sset (\Om^{\,}\a)_t$ extends, via formula
\eqref{f^-}, to a graded derivation $d_t: (\Om A)_t\to(\Om A)_t$
of degree zero.

A useful alternate point of view on the  extended
de Rham complex is provided by the following

\begin{lem}\label{alternate} \vi There is a canonical
DG algebra isomorphism 
$$\xymatrix{
\Om_t A \ar@{=}[r]^<>(0.5){\psi}& (\Om A)_t,
}
$$
such that the de Rham differential on $\Om_t A=\Om_{\k[t]}(A_t)$
is transported, via the isomorphism $\psi$,
to the map $d_t$, see Notation \ref{f_tnot}, on $(\Om^{\,}A)_t$.

\vii The composite isomorphism below
\begin{align*}
\DR A\oplus \Om^{\,}\a
\oplus\! &\xymatrix{{(\Om^{\,}\a)^{\o 2}\cyclic \oplus
(\Om^{\,}\a)^{\o 3}\cyclic \oplus\ldots\;}
\ar@{=}[rr]^<>(0.5){\text{\em Lemma \ref{triv}}}&&
\DR_tA
}
\\
&\xymatrix{\ar@{=}[rr]^<>(0.5){\R(\psi)}&&\R((\Om^{\,}A)_t)
\ar@{=}[rr]^<>(0.5){\text{\em cf. \eqref{dsd}}}&&\oplus_{k\geq 0}
\R^{2k}((\Om A)_t)
}
\end{align*}
is compatible with  the  leftmost and
 the rightmost  direct sum decompositions.
\end{lem}
\begin{proof} In general, for any two algebras $A,R$, we have
a natural DG algebra isomorphism $\Om(A*R)\cong (\Om A)*(\Om R).$
Since $\Om_R(A*R)$ is a quotient of the algebra
$\Om(A*R)$ by the two-sided ideal generated by 
the space $dR\sset \Om^1R\sset\Om^1(A*R)$,
the isomorphism above induces  a  DG algebra isomorphism
$$\Om_R(A*R)\cong (\Om A)*R.$$

Part (i) of the Lemma  follows from the last isomorphism
by putting $R=\k[t]$.
Part (ii) is now an immediate consequence of Lemma
\ref{triv} applied to the algebra $B=\Om A$.
\end{proof}

Yet another interpretation of the extended Karoubi-de Rham complex
will be given in \S\ref{another}.

Next, we observe  that
the algebra $\Om_t A=\Om_{\k[t]}(A_t)$  comes equipped with  a natural 
{\em bi}-grading
$\Om_t A=\oplus_{p,q\geq 0}\, \Om^{2p,q}_t A$, where
the even $p$-grading is induced
from the one on $A_t$, and the $q$-component
corresponds to the grading induced by the
natural one on $\Om^\hdot A$. It is easy to see that
the  $p$-grading corresponds,
under the isomorphism of Lemma \ref{alternate},
to the grading on $(\Om A)*\k[t]$ that
counts twice the number
of occurrences of the variable $t$, and the $q$-component
corresponds to the grading induced by the
natural one on $\Om^\hdot A$,
e.g. for any $\al\in \Om^k A,\be\in \Om^\ell A$ 
the element $\al\,t\,\be\,t \in (\Om A)*\k[t]$
has bi-degree $(2p=4,q=k+\ell)$.

The bi-grading on $\Om_t A$ clearly descends to
a similar bigrading on the
extended de Rham complex of $A,$
with the de Rham  differential of  bi-degree~$(0,1)$:
\beq{grading}\DR_t A=\oplus_{p,q}\, \DR^{2p,q}_t A,\quad
d:\
 \DR_t^{2p,q} A\to\DR^{2p,q+1}_t A.
\eeq

\section{Noncommutative calculus}
\subsection{Lie derivative and contraction for derivations.}
Fix an algebra $B$ and a subalgebra $R\sset B$.

Any  derivation $\th\in \Der_RB$
 gives rise to a  Lie derivative
map $L_\th:  \Om^\hdot_RB$
$\to \Om^\hdot_RB,$
and also to contraction (with $\th$) map
$i_\th: \Om^\hdot_RB\to \Om^{\idot-1}_RB.$
The map $L_\th$  is a degree zero derivation of
the graded algebra $\Om^\hdot_RB$. It is defined
on 1-forms by the formula $L_\th(u\,dv)=(\th(u))\,dv+ u\,d(\th(v)),$
 and then extended uniquely 
 to a map
$L_\th: \Om^\hdot_RB\to \Om^\hdot_RB$
as a derivation.
The  contraction map $i_\th$ is a degree 
$-1$ graded derivation. It is
defined on 1-forms by  formula
 \eqref{der_def} and is
extended to a map
$\Om^\hdot_R B \to\Om^{\idot-1}_R B$
as a
{\em graded}-derivation.

One has the Cartan formula $L_\th=i_\th\ccirc d+ d\ccirc i_\th,$
for any  $\th\in \Der_RB$.
The maps $L_\th$ and $i_\th$ both descend to
well-defined operations on the de Rham complex
$\DR^\hdot_R B=\R^\hdot(\Om_RB),$
satisfying a similar Cartan formula.

Now, let $A$ be an algebra, $\th\in\Der A$ a derivation,
and $L_\th: \Om^\hdot A\to \Om^\hdot A$ the corresponding Lie derivative.

On one hand, the derivation $L_\th$ gives rise, according
to Lemma \ref{f_tlem}(ii),
to a  derivation $(L_\th)_t:
(\Om^{\,}A)_t\to (\Om^{\,}A)_t$.
On the other hand, one may first extend $\th$ to a derivation
$\th_t\in \Der_{\k[t]}(A_t),$ and then consider the corresponding
  Lie derivative
$$L_{\th_t}:\
\Om_{\k[t]}(A_t)=\Om_tA\too\Om_{\k[t]}(A_t)=\Om_tA.$$
Very similarly, we may form
 graded derivations $(i_\th)_t$ and
$i_{\th_t}.$

An important point is that the 
two procedures just explained agree with each other in the
sense that: {\em Under the identification
$\Om_{\k[t]}(A_t)\cong(\Om^{\,}A)_t,$
provided by the isomorphism $\psi$ of Lemma \ref{alternate}(i),
 we have}
\beq{twoconstr}
L_{\th_t}=(L_\th)_t,\quad
i_{\th_t}=(i_\th)_t.
\eeq

All the above operations on $\Om_tA$ descend to
 the extended de Rham complex
$\DR_tA.$
The Cartan identity  $L_\th=i_\th\ccirc d+ d\ccirc i_\th$
induces   similar
identities for the corresponding operations on
$\Om_tA$ and on $\DR_tA$.

\subsection{Lie derivative and contraction for double derivations.}
\label{li_sec}
Associated with any  double derivation  $\Th\in\dder\a$, one has an
$A$-bimodule map 
\beq{iom1}
i_\Th:\
 \Om^1\a\to A\o A,\quad \al\mto i_\Th\al=(i_\Th'\al)\o
(i_\Th''\al),
\eeq
that corresponds to $\Th$ under the  canonical bijection \eqref{der_def}.
With this understood, we have

\begin{prop}\label{chprop} 
Any double derivation $\Th\in\dder\a$ gives rise
canonically to
 a  degree $0$
double derivation
$\BL_\Th\in \dder^0(\Om^\hdot \a),$ and also to a degree $-1$
{\em graded} double derivation $\bbi_\Th\in\dder^{-1}(\Om^\hdot  \a)$.

Explicitly, the  double derivation $\BL_\Th$ is given,
for any $a_0,a_1,\ldots,a_n\in A,$ by the  formula
\begin{align*}
a_0\,da_1\ldots \,da_n&\mto \BL_\Th(a_0\,da_1\ldots \,da_n):=
(\Th'a_0)\o  (\Th''a_0)\,da_1\ldots da_n\\
&+\sum_{1\leq k \leq n}\Bigl(a_0\,da_1\ldots
da_{k-1}\,d(\Th'a_k)\o  
(\Th''a_k)
\,da_{k+1}\ldots da_n\Bigr.\\
&\en +\Bigl.
a_0\,da_1\ldots
da_{k-1}\,(\Th'a_k)\o  
d(\Th''a_k)
\,da_{k+1}\ldots da_n\Bigr).
\end{align*}

The {\em graded} double derivation $\bbi_\Th$ is given,
  for any
 $\al_1,\ldots,\al_n\in \Om^1\a,$  by
 the  formula, cf. \eqref{iom1}:
\begin{align*}
&\al_1\al_2\ldots\al_n\mto
\sum_{1\leq k \leq n}(-1)^{k-1}\cd
\al_1\ldots\al_{k-1}\,(i_\Th'\al_k)\o (i_\Th''\al_k)\,
\al_{k+1}\ldots\al_n.
\end{align*}
\end{prop}

Thus, we have
$$
\BL_\Th:
\Om^n\a\to\bigoplus_{0\leq k \leq n}
\Om^{k}\a\otimes\Om^{n-k}\a,\en\;
\bbi_\Th:
\Om^n \a\to\bigoplus_{1\leq k \leq n}\,
\Om^{k-1} \a\otimes \Om^{n-k} \a.
$$

It is not too difficult to verify by direct computation that the explicit
formula for $\BL_\Th$, resp. for $\bbi_\Th,$ given in
Proposition~\ref{chprop}
does define a degree zero
double derivation $\Om^{\,}\a\to\Om^{\,}\a\o\Om^{\,}\a$, resp. a  degree $-1$
graded double derivation. It is, however, much more
rewarding to  use a more conceptual approach to the
construction of $\BL_\Th$ and $\bbi_\Th$, explained below.

A key point of our approach is that although
double derivations do {\em not}
give rise to  natural operations on the DG algebra
$\Om^\hdot \a$  they do give rise to  
canonical Lie derivative  and contraction
 operations on the {\em extended} DG algebra $\Om^\hdot_tA$.
Proposition \ref{chprop} will easily follow from this.

To construct  Lie derivative  and contraction
 operations associated with
 a double derivation, we are going to mimic
equations \eqref{twoconstr}. Observe that the left-hand sides
of equations \eqref{twoconstr} still make sense for double derivations.
In more detail, given  $\Th\in \dder{A}$,  we may first extend it
to a free product derivation $\Th_t: A_t\to A_t$,
as in \S\ref{free_prod}. After that,
we define  Lie derivative $L_\Th$, resp. contraction
$i_\Th$,  operation on
$\Om_t A$, as the ordinary Lie derivative, resp.  ordinary contraction,
operation induced  by the derivation $\Th_t$
on the complex $\Om_t A=\Om_{\k[t]}(A_t),$
of {\em relative} differential forms on the algebra $A_t$.
Specifically, we put
\beq{li}
L_\Th:=L_{\Th_t}\en(\text{resp.},\en
i_\Th:=i_{\Th_t})\, :\
\Om_t A\map \Om_t A.
\eeq

Now, we use Lemma \ref{triv} to write $\Om_tA=(\Om^{\,}A)_t$.
It is immediate from the above formulas  that,
viewed as  maps  $(\Om^{\,}\a)_t\to (\Om^{\,}\a)_t,$ 
the Lie derivative $L_\Th$ and  contraction $i_\Th$ are both
graded derivations of degree 2 with respect to the 
grading that counts
 (twice) the number of occurences of $t$
(and disregards
the degrees of differential forms).
Therefore, applying Lemma \ref{gradings} to the
algebra $B=\Om\a$, we deduce that
there exist   unique double
derivations $\BL_\Th,\bbi_\Th:  \Om^{\,}\a\to (\Om^{\,}\a)\o (\Om^{\,}\a),$
such that,
for the corresponding maps $(\Om^{\,}\a)_t\to(\Om^{\,}\a)_t,$ we have,
cf.  \eqref{twoconstr}
\beq{L^T}
(L_\Th:=)\,L_{\Th_t}=(\BL_\Th)_t,\quad\text{resp.}\quad 
(i_\Th:=)\,i_{\Th_t}=(\bbi_\Th)_t.
\eeq
(in the case of contraction operation $i_\Th$, one has to repeat the argument
using an analogue of Lemma \ref{gradings} for {\em graded}-derivations).

This completes our proof of Proposition
\ref{chprop}. \qed\smallskip

The Lie derivative 
 and contraction operations on $\Om_tA$ defined in \eqref{li}  satisfy,
for any $\Th,\Phi\in\dder\a$ and $\xi\in\Der A$, the following standard
commutation relations
\beq{relations}
L_\Th=d\ccirc i_\Th+i_\Th\ccirc d,
\quad
i_\Th\ccirc i_\Phi+ i_\Phi\ccirc i_\Th=0,
\quad i_\xi\ccirc i_\Th+i_\Th\ccirc i_\xi=0.
\eeq

To prove \eqref{relations}, one first verifies these identities
on the generators of the algebra $\Om_t A=(\Om^{\,}\a)_t$,
that is, on differential forms of degrees $0$ and $1$,
which is a simple computation.
The general case then follows by observing
that any commutation relation between (graded)-derivations
that holds on generators of the algebra holds true on
all elements of the algebra.

Formulas \eqref{relations} imply, in particular, that the Lie derivative $L_\Th$
commutes with the de Rham differential $d$.

Both the Lie derivative and contraction operations on $\Om_tA$ 
descend to the commutator quotient
$\DR_t A$. This way, we obtain the Lie derivative
 $L_\Th$, resp.  contraction
$i_\Th,$  operation on the extended de Rham complex.
It is immediate that these operations 
on $\DR_t A$ also satisfy \eqref{relations}.

\subsection{Reduced Lie derivative and contraction.}\label{red_sec}
Let $B$ be a (possibly graded)
 algebra and $F\in\Der_t^2(B_t)$  a degree 2 derivation
with respect to the {\em even} grading on $B_t$ that counts twice 
the number of occurencies of $t$ and disregards the grading on $B$.
Using the identifications of Lemma \ref{triv}, we can interpret
the induced linear map  $\R(F) : \R^\hdot(B_t)\to\R^{\idot+2}(B_t)$ as
a chain of maps
\beq{T_t}
B/[B,B]\,\stackrel{\R(F)}\too\,
B\,\stackrel{\R(F)}\too\, B^{\o 2}\cyclic\,\stackrel{\R(F)}\too\,
B^{\o 3}\cyclic \,\stackrel{\R(F)}\too\,\ldots\,.
\eeq

We apply the above in the special case where  $B=\Om^{\,}\a$ and
 $F:=\BL_\Th$, resp. $F:=\bbi_\Th$, for some $\Th\in\dder\a$.
In this case, the terms in \eqref{T_t} may be written
more explicitly using  the identification provided by the
composite isomorphism of  Lemma
\ref{alternate}(ii).
This way, any double derivation $\Th\in\dder\a$ gives rise
to two chains  of maps (one for $\BL_\Th$ and one for $\bbi_\Th$)
of the form
\beq{chainDR}
\DR A\too \Om^{\,}\a\too   
(\Om^{\,}\a)^{\o 2}\cyclic \too 
(\Om^{\,}\a)^{\o 3}\cyclic \too \ldots.
\eeq

Recall now that the extended de Rham complex $\DR_tA$ 
has a natural bigrading, cf.~\eqref{grading}. The second component
(the $q$-component) of
that
bigrading induces
a grading
on each of the  spaces $\dis(\Om^{\,}\a)^{\o k}\cyclic,\,k=1,2,\ldots,$
appearing in \eqref{chainDR}.
Explicitly,
the grading on $\dis(\Om^{\,}\a)^{\o k}\cyclic$
counts the total degree of differential
forms involved, i.e., we have $\dis\deg(\al_1\o\ldots\al_k)=
\deg\al_1+\ldots+\deg\al_k,$ for any $\dis\al_1\o\ldots\al_k
\in (\Om^{\,}\a)^{\o k}\cyclic.$
Further, it is immediate from  definitions that
the maps $L_\Th$ and $i_\Th$ in \eqref{li} are graded derivations
of bi-degrees
$(2,0)$ and $(2,-1),$  respectively. 
We conclude that, in the  Lie derivative case,
all maps in the corresponding chain  \eqref{chainDR}
preserve the above defined grading
while, in the contraction case,
all maps in the corresponding chain  \eqref{chainDR}
decrease  the  above defined grading by one.

The  leftmost map in \eqref{chainDR}, to be denoted $\bi_\Th$
in the contraction case, resp.
$\bl_\Th$ in the Lie derivative case,
will be especially important for
us. In more detail,  these maps have the form
\beq{i_explicit}
\bi_\Th: \ \DR^\hdot A\map \Om^{\idot-1}A,
\quad\text{resp.},\quad
\bl_\Th: \
\DR^\hdot A\map \omb.
\eeq
One verifies that the map $\bi_\Th$ in \eqref{i_explicit} is given,
for any $\al_1,\al_2,\ldots,\al_n\in\Om^1 A,$ 
by the following explicit formula, see \cite[(2.8.4)]{CBEG}:
\begin{align}\label{bi_formula}
\bi_\Th(\al_1&\al_2\ldots\al_n)=\\
&=
\sum\nolimits_{k=1}^n
(-1)^{(k-1)(n-k+1)}\cd  (i_\Th''\al_k)\cd
\al_{k+1}\ldots\al_n\,\al_1\ldots\al_{k-1}\cd (i_\Th'\al_k).\nonumber
\end{align}

The maps in \eqref{i_explicit}  have been already introduced,
in effect, 
in \cite{CBEG} under the names of {\em reduced contraction},
resp. {\em reduced Lie derivative}.
An {\em ad hoc}  definition of these maps given in \cite{CBEG} was based on
 explicit formulas for $L_\Th$ and $i_\Th$, cf. formulas
 of Proposition \ref{chprop}.
As a result, verification of commutation relations
like \eqref{relations} was very painful; it was carried out
in  \cite{CBEG} by rather long brute force computations.
Our present approach based on the free product construction
yields the  commutation relations almost for free.

\subsection{The derivation $\Delta$.}\label{vdb}
For any algebra $A$, there is a distinguished
derivation 
$$\Delta: A\to A\otimes A,\en
a\mapsto 1\o a-a\o 1.
$$

The corresponding contraction map
$i_\Delta: \oma\to  A\otimes A$ is nothing but the tautological
imbedding $\oma\into A\otimes A$. Further, the derivation
$\Delta_t: A_t\to A_t$ associated with $\Delta$ by formula \eqref{f_t} clearly
equals $\ad t: u\mapsto t\,u-u\,t$.
Given $a\in A$ and $\om\in\Om^{\,}\a,$  write $[\om,a]:=\om\,a-a\,\om.$ 

\begin{lem}\label{biDelta} \vi The Lie derivative
$L_\Delta: \Om_tA\to\Om_tA$ is given by the formula
$\om\mapsto [t,\om]$.
In particular, 
the induced map $\R(L_\Delta):  \DR_tA\to\DR_tA,$ on the extended de Rham
complex, {\em vanishes}.

\vii For any $a_0,a_1,\ldots, a_n\in A,$ we have
\begin{align*}
\bi_\Delta(a_0\,da_1\ldots \,da_n)
=\sum_{1\leq k \leq n}(-1)^k\cd
[da_{k+1}\ldots da_n\,a_0\,da_1\ldots da_{k-1},\;a_k].
\end{align*}
\end{lem}
\begin{proof}  To prove (i),
it suffices to check the equality $L_\De=\ad t$ on the generators
of the algebra $\Om_tA$, that is, on 0-forms and on 1-forms.
This is easy. Part (ii) is a straightforward consequence of 
formula \eqref{i_explicit}.
\end{proof}

There is also an {\em odd} version of the derivation $\Delta_t$.
Specifically, one can extend the double derivation $\Delta$ to
a degree $-1$ graded derivation
$\Delta_\tau: B_\tau\to B_\tau$ (cf. \S\ref{ttau} for the definition of
$B_\tau$), such that
$\Delta_\tau(b)=\tau\,b-b\,\tau, \,\forall b\in B\sset
B_\tau,$
and $\Delta_\tau(\tau)= 0$. 
Further, it is straightforward to check that
$$ \mbox{$\frac{d}{d\tau}$}\ccirc\Delta_\tau
+\Delta_\tau\ccirc\mbox{$\frac{d}{d\tau}$}=0,
\quad\text{and}\quad \Delta_\tau\ccirc\Delta_\tau=-\ad(\tau^2),
$$
where $\ad(\tau^2)$ stands for the graded commutator with $\tau^2$.
Observe that $\Delta_\tau\neq \ad\tau,$ since
 $\Delta_\tau(\tau)=0\neq 2\tau^2=\ad\tau(\tau).$

We see that the pair $\frac{d}{d\tau}$ and $\Delta_\tau$
gives rise, by Proposition \ref{drin}(ii), to two anti-commuting
maps of degrees $+1$ and $-1$, respectively,
on the space $\R^{<0}(A_\tau)$, hence, on the 
cyclic complex \eqref{cyc_comp}.
 Using the isomorphism of  Proposition \ref{drin}
we find that 
the corresponding map $A\mcyclic^{\o n}\to
A\mcyclic^{\o (n+1)},$ arising from $\frac{1}{2}\R(\Delta_\tau),$
is given by the following explicit formula
$$\mbox{$\frac{1}{2}$}\R(\Delta_\tau)(a_1\o\ldots\o a_n)=
\sum\nolimits_{k=1}^n (-1)^ka_1\o\ldots\o a_k\o1\o a_{k+1}\o\ldots\o a_n.
$$

\section{The Representation functor.}\label{repf}
\subsection{Evaluation map.}\label{ev_sec}
 We  fix  a finite dimensional
$\k$-vector space $V$, and write
$\End:=\Hom_\k(V,V)$ for the associative algebra of
linear endomorphisms of $V$.

Given an algebra $A$, one may consider
the set $\Hom_{\mathsf{alg}}(A,\End)$ of all algebra maps
$\rho: A\to \End$.
More precisely, to any {\em finitely generated}
$\k$-algebra $A$ one associates
an affine scheme  of finite type
over $\k$,
to be denoted $\Rep(A,V),$
such that the set  $\Hom_{\mathsf{alg}}(A,\End)$
is the set of $\k$-points of  $\Rep(A,V)$.
Write $\k[\Rep (A,V)]$ for the
coordinate ring of the affine scheme $\Rep (A,V),$
which will be always assumed to be  {\em non-empty}.

The  tensor product  $\End\o\k[\Rep (A,V)]$,
is an associative
algebra
that
may also be identified with the algebra of polynomial
maps $\rep\to\End$, equipped with pointwise multiplication.

To each element $a\in A$, one associates
the  function $\wh{a}: \rep\to\End,$
$\rho\mapsto \wh{a}(\rho):=\rho(a)$. The assignment
$a\mapsto\wh{a}$ clearly gives an algebra
homomorphism, called  {\em evaluation map}, cf. also \cite[\S12]{Gi},
$$ \ev:\ A\too\End\o\k[\Rep (A,V)],
\quad a\mto\wh{a}.
$$

\subsection{Extended de Rham complex and equivariant cohomology.}
In this section, we assume  the field $\k$ to be algebraically closed.

Let $X$ be an affine scheme with coordinate ring
$\k[X]$. 
We write
$\calT(X):=\Der(\k[X])$ for
the  Lie algebra of algebraic vector fields,
$\Om^1(X)$ for the $\k[X]$-module of
K\"ahler differentials, and let
 $\Omega^\hdot(X):=\bigwedge^\hdot_{\k[X]}\Om^1(X)$ 
denote the DG algebra of differential forms on $X$.

Suppose we have a finite dimensional Lie algebra $\g g$ and
a Lie algebra  map $\g g\to\calT(X)$. 
Thus, we get a Lie algebra action of $\g g$
on $\Omega^\hdot(X),$ by the Lie derivative.
Also, let $\g g$ act on $\C[\g g]$,
the polynomial algebra on the vector space $\g g$, by the
adjoint action. We view  $\C[\g g]$ as an even-graded
algebra such that the vector space of linear functions
on $\g g$ is assigned degree 2.
This makes the
tensor product $\Omega^\hdot(X,\,\g g):=\Omega^\hdot(X)\otimes\C[\g g]$
a graded algebra, equipped with
$\g g$-diagonal action. We write
$\Omega^\hdot(X,\,\g g)^{\g g}$ for the corresponding  graded subalgebra
of $\g g$-invariants.

The $\g g$-{\em equivariant algebraic de Rham complex} of $X$ is  the complex
\beq{dr_equiv}
\bigl(\Omega^\hdot(X,\,\g g)^{\g g},\,d_{\text{DR}}+d_\g g).
\end{equation}
Here,
the differential  $d_{\text{DR}}=d\o \id_{\C[\g g]}$ is induced by the
ordinary de Rham 
differential $d$ on $\Omega^\hdot(X)$. The differential  
$d_\g g$ 
 is defined by the formula
$$
d_\g g:\
 \om\o f\mto \sum\nolimits_{r=1}^{\dim \g g} (i_{e_r}\om)\o (e_r^*\cd f),
\quad\forall \om\o f\in \Omega^\hdot(X)\otimes\C[\g g],
$$
where  $\{e_r\}$ and  $\{e_r^*\}$ stand for  dual bases  of $\g g$ and $\g g^*$,
respectively,
and
where  $i_{e_r}$ denotes contraction by the vector field on $X$
corresponding
to the element $e_r$.

We now return to the setup of \S\ref{ev_sec}. Thus we fix a finitely
generated
algebra $A$, a finite dimensional vector space $V$,
and consider the scheme $\Rep(A,V).$

The  general linear group $\GL(V)$ acts naturally on
the algebra $\End=\Hom_\k(V,V)$ by inner automorphisms, via
conjugation. Hence, given an algebra homomorphism
$\rho: A\to \End$ and $g\in \GL(V)$, one
may define a conjugate homomorphism
$g(\rho): a\mapsto g\cdot\rho(a)\cdot g\inv.$
The action $\rho\mapsto g(\rho)$ makes 
$\rep$ a $\GL(V)$-scheme. 

Let $\g g:=\Lie \GL(V)$ be the Lie algebra of
 the algebraic group $\GL(V)$.
The action of  $\GL(V)$ on $\rep$ induces
a Lie algebra  map 
\beq{actx}
\act_A:\ \g g \map\calT(\Rep(A,V)),
\quad x\mto \act_A(x).
\eeq
Thus, one may consider $\Omega^\hdot(\Rep(A,V),\,\g g)^{\g g},$
 the corresponding 
 $\g g$-equivariant algebraic de Rham complex.

Now,  we have defined a contraction operation
$i_\Delta: \DR_t A\to\DR_t A,$ cf. \S\ref{li_sec}.
The main result of this section is the following
\begin{thm}\label{rep_thm}
There is a canonical morphism of complexes compatibale with both differentials
$$\bigl(\DR^\hdot_t A,\, d,\,i_\Delta\bigr)\stackrel{\ev_t}\too
\bigl(\Omega^\hdot(\Rep(A,V),\,\g g)^{\g g},\,d_{\text{DR}},\,d_\g
g\bigr).
$$
\end{thm}

We begin the proof with some general constructions. 

\subsection{Evaluation map on differential forms.}
Let $B$ be an algebra and
$R\sset B$ a subalgebra. The algebra imbedding $R\into B$ induces,
by transport of structure,
a canonical morphism of schemes $\Rep(B,V)\to\Rep(R,V).$

Observe next that giving an algebra homomorphism
$\rho: \k[t]\to \End$ amounts to  specifying an
element $x=\rho(t)\in\End,$ that may
be chosen arbitrarily. Thus, we have $\Rep(\k[t],V)=\End$.

Similarly, for any algebra $A$, giving an algebra morphism
$\rho:A*\k[t]$
$\to\End$ amounts to giving a 
homomorphism
$A\to \End$ and an arbitrary additional element $x=\rho(t)\in\End.$
We see that
$\Rep(A_t,V)=\Rep(A,V)\times\End.$
With this identification, the canonical morphism $\Rep(A_t,V)\to \Rep(\k[t],V),$
induced by the  algebra imbedding $\k[t]\into A_t,$
is nothing but the second projection
$\pi: \Rep(A,V)\times\End\to\End.$

We put $G:=\GL(V).$ This group acts naturally  on each of the schemes
$\Rep(A_t,V),\,\Rep(A,V),\,\Rep(\k[t],V),$
and the projection $\pi$ is clearly $G$-equi-variant.
We will use shorthand notation 
$$\Rep_t:=\Rep(A_t,V)=\Rep(A,V)\times\End,\quad\text{\em
and}\quad\Rep:=\Rep(A,V).
$$

Let
$\Om^\hdot_\pi(\Rep_t)$
be the DG algebra of relative (with respect to $\pi$)
 algebraic differential forms on
the scheme $\Rep_t$ (in the ordinary sense of commutative algebraic
geometry).
By definition, we have
\beq{rel}
\Om^\hdot_\pi(\Rep_t):=\wedge^\hdot_{\k[\Rep_t]}\,\Om^1_\pi(\Rep_t)=
\Om^\hdot(\Rep)
\otimes \k[\End].
\eeq 

Generalizing the construction of \S\ref{ev_sec},
 we now construct an evaluation map
$\dis
\Om^\hdot_{\k[t]}(A_t) \too\End\o\Om^\hdot_\pi(\Rep_t),$
 on relative
differential forms as follows. 

First of all, given $u\in A_t$, we may differentiate
the $\End$-valued function $\wh{u}: \Rep_t\to\End$
to obtain an  $\End$-valued 1-form $d\wh{u},$
on $\Rep_t.$ The assignment
$u\mapsto d\wh{u}$ gives a map
$A_t\to \End\o\Om^1(\Rep_t),$ which
descends to a well-defined
map $d\ccirc \ev:
A_t/\k[t]\to \End\o\Om^1_\pi(\Rep_t)$.

Now, given $n=0,1,\ldots,$
write $m:$
 $\End^{\otimes (n+1)}\to\End$ for the $(n+1)$-fold
multiplication map in the algebra $\End$
and  ${\text{Alt}: T^n(-)\to\wedge^n(-)}$ for
the projection to anti-symmetric tensors.
We form the following composite
\begin{align*}
&\xymatrix{\Om^n_{\k[t]}(A_t)\ar@{=}[r]& 
{A_t\otimes_{\k[t]} T^n_{\k[t]}(A_t/\k[t])}
\ar[rrr]^<>(0.5){\ev\o (d\ccirc\ev)^{\o n}}&&&}\\
&(\End\o\k[\Rep_t])\otimes_{\k[t]}
T^n_{\k[t]}(\End\o\Om^1_\pi(\Rep_t))=\\
&{\xymatrix{
\End^{\otimes n+1}\o
T^n_{\k[t]}(\Om^1_\pi(\Rep_t))
\ar[rr]^<>(0.5){m\o\text{Alt}}&&
\End\o
\bigl(\wedge_{\k[\Rep_t]}^n\,\Om^1_\pi(\Rep_t)\bigr)
}}\\
&\qquad\qquad\qquad\qquad\qquad\qquad\qquad\qquad
=\End\o\Om^n_\pi(\Rep_t).
\end{align*}

All these maps  clearly commute with the de Rham differentials.
Furthermore, any element in the image of the composite 
above is easily seen to
be $G$-invariant with respect to the $G$-diagonal action on
$\End\o\Om^n_\pi(\Rep_t).$ Thus, the composite  yields
 a well-defined, canonical DG  algebra map
$$\evom:
\Om_t A\to
\bigl(\End\o\Om^\hdot_\pi(\Rep_t)\bigr)^G\!,
\en\;
\al=a_0\,da_1\ldots da_n\mapsto \wh{\al}=\wh{a}_0\,d\wh{a}_1\ldots d\wh{a}_n.$$

Further, we have the linear function
 $\Tr: \End\to\k, \,x\mapsto\Tr(x).$
We  form the  following composite
$$\Om_t A\stackrel{\evom}\too
\bigl(\End\o\Om^\hdot_\pi(\Rep_t)\bigr)^G\stackrel{\Tr\o\Id}\tooo
 \bigl(\k\o\Om^\hdot_\pi(\Rep_t)\bigr)^G=\Om^\hdot_\pi(\Rep_t)^G.$$ 
The map  $(\Id\otimes\Tr)\ccirc\evom$ vanishes on 
$[\Om_t A,\Om_tA]$ $\sset\Om_t A$, the (graded) commutator space,
due to symmetry of the trace function.
Therefore, this map
descends to the commutator quotient to give
  a well-defined $\k$-linear map of graded vector spaces
\beq{phi}(\Id\otimes\Tr)\ccirc\evom:\
\R^\hdot(\Om_tA)\map \Om^\hdot_\pi(\Rep_t)^G,
\quad \al\mapsto
\Tr\wh{\al}.
\eeq
By construction, this map   commutes with the de Rham differentials.

We remark next that the Lie algebra $\g g=\Lie \GL(V)$ is nothing
but the associative algebra $\End$ viewed as a Lie algebra.
Hence, 
using the isomorphisms in \eqref{rel},
we can write 
$$\Om^\hdot_\pi(\Rep_t)=\Om^\hdot(\Rep)\o\k[\End]=\Om^\hdot(\Rep)
\o\k[\g g]=\Om^\hdot(\Rep,\,\g g).$$
Thus, by the definition of the extended de Rham complex, $\DR^\hdot_tA,$
the map in \eqref{phi} may be viewed as a map
\beq{trace_om2}
\ev_t:=(\Id\otimes\Tr)\ccirc\evom:\
\DR^\hdot_tA\map \Om^\hdot_\pi(\Rep_t)^G=
\Om^\hdot(\Rep,\,\g g)^{\g g}.
\eeq

\subsection{Proof of Theorem \ref{rep_thm}.} We are going to show that
the map
\eqref{trace_om2} satisfies the requirements for the
map $\ev_t$ in Theorem \ref{rep_thm}.
Thus, proving the theorem amounts to  showing commutativity
of the following diagram
\beq{diag}
\xymatrix{
\DR_tA\ar[d]_<>(0.5){i_\Delta}\ar[rrr]^<>(0.5){(\Id\otimes\Tr)\ccirc\evom}&&&
\Om(\Rep,\,\g g)^{\g g}\ar[d]_<>(0.5){d_{\g g}}\\
\DR_tA\ar[rrr]^<>(0.5){(\Id\otimes\Tr)\ccirc\evom}&&&
\Om(\Rep,\,\g g)^{\g g}.
}
\eeq

To this end, we need first to recall some generalities
about the effect of 
the Representation functor on derivations,
see e.g. \cite[\S12.4]{Gi} for details.

Let $B$ be an arbitrary finitely generated algebra.
Each derivation
$\delta\in \Der B$ gives rise, in a canonical way,
to a vector field on the scheme $\Rep(B,V)$, that is,
to a derivation $\wh{\delta}: \k[\Rep(B,V)]\to\k[\Rep(B,V)].$
The assignment
$\delta\mapsto\wh{\delta}$ is a Lie algebra 
morphism $\Der B\to \Der\k[\Rep(B,V)].$
This  morphism is compatible with
the Representation functor in the sense that one has
$$\wh{\delta}(\Tr\ccirc\ev(b))=\Tr\ccirc\ev(\delta(b)),
\quad\forall b\in B,\,\delta\in\Der B.
$$

Let $\rho\in\Rep(B,V)$ be a (closed) point. We
write $\wh{\delta}|_\rho\in T_\rho\Rep(B,V)$ for
the value 
of the vector field $\wh{\delta}$ at the point $\rho$;
thus,  $\wh{\delta}|_\rho$ is
an element of the Zariski tangent space to the scheme $\Rep(B,V)$ at $\rho$.
Similarly, for any $x\in\g g=\Lie \GL(V),$ 
let $\act_B(x)|_\rho$ denote the
value at  $\rho$ of the  vector field $\act_B(x)$, see \eqref{actx}.

Now, let $b\in B$ and set $\delta=\ad b$, an {\em inner} derivation.
Then,
according to \cite[Example 12.4.1]{Gi}, 
we have 
\beq{act}
\wh{\ad b}|_\rho= \act_B(\rho(b))|_\rho,\quad\forall\rho\in\Rep(B,V),
\eeq
where $\rho(b)\in\End$ is viewed as an element of $\g g=\Lie
\GL(V)=\End.$

\begin{proof}[{Proof of commutativity of diagram \eqref{diag}.}]
{We will use the identification}\break
${\Rep(A_t,V)=\Rep(A,V)\times\End=\Rep\times\g g}$ to write
an element of $\Rep(A_t,V)$ as a pair $(\rho,x)\in \Rep\times\g g.$
The action of $\g g$ on $\Rep(A_t,V)$ corresponds
to the $\g g$-diagonal action on $\Rep\times\g g$,
where the Lie algebra $\g g$ acts on the second  factor, $\g g$, via
 the adjoint action. Thus, for $x\in \g g$, the corresponding vector
field on $\Rep(A_t,V)$ splits into a sum 
 $\act_{A_t}(x)=\act_A(x)\oplus
\ad_{\g g} x$, of the vector field $\act_A(x)$ on $\Rep$ and of
the adjoint vector field $\ad_{\g g} x$ on $\g g$ defined by the formula
 $\ad_{\g g} x|_y=[x,y],$ for any $y\in\g g$.

We apply formula \eqref{act}  to the algebra
$B=A_t$ and the element $b=t$. Recall that  according to  \S\ref{vdb},
for the derivation $\Delta_t: A_t\to A_t,$ associated with the
distinguished double derivation $\Delta: A\to A\o A$,
we have $\Delta_t=\ad t$. Thus, for any point $(\rho,x)\in \Rep\times\g g,$
formula  \eqref{act} yields
\begin{align}\label{t}
\wh{\Delta}_t|_{(\rho,x)}&=\wh{\ad t}|_{(\rho,x)}=
\act_{A_t}(x)|_{(\rho,x)}\\
&=
(\act_A(x)|_{\rho})\,\oplus\,
(\ad_{\g g} x|_{x})
=(\act_A(x)|_{\rho})\,\oplus\, 0,\nonumber
\end{align}
where in the last equality we have used that 
$\ad_{\g g} x|_x=[x,x]=0.$

We observe that the vector field on the right of the second line
in \eqref{t} is
{\em vertical}, i.e.,  is tangent to the fibers of the
second projection $\pi: \Rep\times\g g\to\g g$.
For any vertical vector field $\xi$ on $\Rep_t$ and
any differential
form $\om\in\Om^\hdot(\g g),$
we clearly have $i_\xi(\pi^*\om)=0$.
It follows in particular that, for $\xi=\wh{\Delta}_t,$
the corresponding contraction operation $i_\xi$
descends to a well defined operation on the space
$\Om^\hdot_\pi(\Rep_t)$ of {\em relative} differential
forms. 

To complete the proof we must show that
$\dis i_{_{\wh{\Delta}_t}}=d_{\g g}$, as maps $\Om(\Rep,\,\g g)^{\g g}\to
\Om(\Rep,\,\g g)^{\g g}.$
To see this, choose dual bases $\{e_r\}$ and $\{e^r\}$ of
the vector spaces $\g g$ and $\g g^*$, respectively.
Thus, for any $x\in\g g$, one has
$x=\sum_r e^r(x)\cdot e_r$. Hence, we can write
$\act_A(x)=\sum_r e^r(x)\cdot \act_A(e_r)$.
Therefore, for   any differential
form $\om\in\Om^\hdot(\Rep_t),$ using formula \eqref{t} and the
definition
of the differential $d_{\g g}$,
we find
$$
(i_{_{\wh{\Delta}_t}}\om)|_{(\rho,x)}=(i_{\act_A(x)}\om)|_{\rho}=
\sum\nolimits_r e^r(x)\cdot (i_{\act_A(e_r)}\om)|_{\rho}=
(d_{\g g}\om)|_{(\rho,x)}.
$$

This completes the proof of commutativity of diagram \eqref{diag},
hence, the proof of the theorem. 
\end{proof}

\section{Applications to Hochschild and  Cyclic homology}
\subsection{Hochschild homology.}\label{hoch_sect}
Let $A$ be a unital $\k$-algebra.
Given an $A$-bimodule $M$, we
let $H_k(A,M)$ denote the $k$-th  Hochschild
homology group of $A$ with coefficients in $M$.
Also, write $[A,M]\sset M$ for the
$\k$-linear span of the set $\{am-ma\mid a\in A,m\in M\}$.
Thus, $[A,M]$ is a vector subspace, and
$H_0(A,M)= M/[A,M]$.

It has been shown by  Cuntz-Quillen  that,
in the special case $M=A,$  the
 Hochschild
homology groups $H_\idot(A,A)$ may be computed
using noncommutative differential forms. Adapting some of the
arguments from  \cite{CQ2}, we will establish our first result.

\begin{thm}\label{hoch_thm} For any unital $\k$-algebra $A$, there
is a natural graded space isomorphism
$$H_\idot(A,A)\cong\Ker[\bi_\De: \DR^\hdot A\to \Om^{\idot-1}A].$$
\end{thm}

As will be explained later (see discussion after  Proposition
\ref{bkappa}), this theorem is an easy consequence of Proposition
\ref{bkappa}; the latter proposition
will be itself proved in Sect. \ref{pf1} below.

To put  Theorem \ref{hoch_thm} in context, following
 Cuntz and Quillen,   consider a complex
$\ldots\stackrel{b}\map\Om^{2}A
\stackrel{b}\map\Om^{1}A\stackrel{b}\map\Om^{0}A\map 0,$
where  the differential $b$ is given by the
formula
\beq{b} b:\ \al\,da\mto(-1)^{n }\cd[\al, a],\quad\forall
 \al\in\Om^{n }A,\,
n >0,\, a\in A/\k
\eeq

It has been explained in \cite{CQ2} that this complex 
is isomorphic naturally to the   Hochschild
complex for the algebra $A$. Therefore,
the cohomology groups $H^{\idot}(\Om^{\,}\a,b),$ of the above complex,
may be identified with   Hochschild
homology groups $H_\idot(A,A)$.

\subsection{Cyclic homology.}
We introduce a complex of {\em reduced} differential forms
 by setting 
$\bom^{_0}:=\Om^0A/\k=A/\k$, and
$\bom^{_k}:=\Om^kA,$ for all $k>0$.
Let $\bom^{_\idot} :=\bigoplus_{k\geq 0}\,\bom^{_k}$.

Further, set $\bom[t]:=\bom^{_\idot}\otimes\k[t]$,
viewed as a graded vector space with respect to the
tensor product grading where the variable $t$ is assigned 
grade degree 2, as usual.

The Hochshild differential induces a $\k[t]$-linear
differential $b: \bom[t]\to\bom[t]$ of degree $(-1)$.
One also has Connes' differential
$B: \bom^{_\idot} \to\bom^{_{\idot+1}}$,
and we  extend it to a $\k[t]$-linear 
differential $\bom[t]\to\bom[t]$ of degree $+1$.
It is known that $B^2=b^2=0$ and 
$B\ccirc b+b\ccirc B=0$. Thus, the map
$B+t\cdot b: \bom[t]\to\bom[t]$
gives a degree $+1$ differential on $\bom[t]$.

The  reduced cyclic homology of $A$,
to be denoted $\HC_\idot(A)$, is defined
as homology of the complex
$(\bom[t], B+t\cdot b)$.

On the other hand, we have the reduced contraction
 $\bi_\Delta: \DR^\hdot A\to\Om^{\idot-1}A$. We compose it
with the natural projection $\omb\onto \DR^\hdot A$
to obtain a map  $\omb\to\Om^{\idot-1}A$.
The latter map descends to a map
 $\bom^{_\idot} \to\bom^{_{\idot-1}}$.
Further, we may extend this last map, resp., de
Rham differential
$d: \bom^{_\idot} \to\bom^{_{\idot+1}},$
to  $\k[t]$-linear maps
 $\bom[t]\to\bom[t]$ of degree $-1$,
resp., $+1$,
to be denoted by $\bi_\Delta,$ resp. $d$, again.
It follows from Lemma \ref{biDelta}(i) and the Cartan formula,
cf. also  \eqref{relations}, that one has
\beq{iddi} 
i_\Delta\ccirc d+d\ccirc i_\Delta=0,\quad
d^2=(\bi_\Delta)^2=0\quad
\text{on}\en \DR_tA.
\eeq
Thus, the map $d+t\cdot\bi_\Delta$
gives a degree $1$ differential on $\bom[t]$ that may be thought of as
some sort of equivariant differential for the `vector field'
$\Delta$.

The following theorem, which is the main result of this paper,
shows the importance of the reduced contraction
map $\bi_\Delta$ for Cyclic homology.

\begin{thm}\label{main} The homology of the  complex
$(\bom^{_\idot}[t], d+t\cdot\bi_\Delta)$ 
is isomorphic to
$\HC_\idot(A)$, the reduced cyclic homology of $A$.
\end{thm}

The proof of the theorem is
 based on various results about
the Karoubi operator  to be  discussed first.
Theorem \ref{main} will be proved after that, in \S\ref{pf2} below.

\subsection{The Karoubi operator.}\label{hoh} 
For any algebra $A$ and an $A$-bimodule $M$, 
 we put
$M\br:= M/[A,M]=H_0(A,M).$ 
Now, let $A\to B$ be  an algebra homomorphism.
Then, $B$ may be viewed as an $A$-bimodule, and 
one has a canonical projection $B\br=B/[A,B]\onto \R B=B/[B,B]$.
In particular, for $B=\omb$, we get
a natural projection $(\omb)\br\to\DR^\hdot A$
which is not 
an isomorphism, in general.

Following Cuntz-Quillen \cite{CQ2},
we consider a diagram
$$
\xymatrix{
{{\Om^{0}}A\;} \ar@<1ex>[r]^<>(0.5){d}&
{\;{\Om^{1}A}\;} \ar@<1ex>[r]^<>(0.5){d}\ar[l]^<>(0.5){b}&
{\;{\Om^{2}A}\;} \ar@<1ex>[r]^<>(0.5){d}\ar[l]^<>(0.5){b}&
\ar[l]^<>(0.5){b}\ldots.
}
$$
where  the   {\em Hochschild differential} $b$
has been defined in \eqref{b}.

The maps $d$ and $b$ are related via an important {\em  Karoubi operator}
 $\ka: {\Om^\hdot  A}\to{\Om^\hdot  A},$ defined by the formula
$\ka:\al\,da\mto(-1)^{\deg\al}\,da\,\al\,$ if $\deg\al>0$,
and $\ka(\al)=\al$ if $\al\in {\Om^0 A}$.
One has, see \cite{CQ1}.
$$b\ccirc d +d\ccirc b=\Id-\ka.$$

It follows that $\ka$ commutes with both $d$ and $b$.
Further, it is easy to verify, cf. \cite{CQ1} and the proof of Lemma
\ref{ka_iso} below, that the Karoubi operator
descends
to a well-defined map $\ka: (\Om^n  A )\br\to(\Om^n  A )\br$,
which is essentially a cyclic permutation;
specifically, we have  
$$\ka(\,\al_1\,\al_2\,\ldots\,\al_{n-1}\,\al_n\,)=
(-1)^{n-1}\al_n\,\al_1\,\al_2\,\ldots\,\al_{n-1},\quad\forall
\al_1,\ldots,\al_n\in\oma.
$$

\begin{prop}\label{bkappa} For any $n\geq 1$, we have
an equality
$$\bi_\Delta=(1+\ka+\ka^2+\ldots+\ka^{n-1})\ccirc b
\quad\op{as}\;\;\op{maps}\en{{\Om^n  A}}\;\to\;{\Om^{n-1}  A}.
$$

Furthermore, the map $\bi_\Delta$ fits into a canonical
short exact sequence
$$0\map H^n(\Om^{\,}\a,b)\map \DR^n A\stackrel{\bi_\Delta}\too 
[A,\Om^{n-1} A]^\ka\map 0.
$$
\end{prop}

Here and below, $(-)^\ka$ denotes taking $\ka$-invariants,
in particular, we write $(\Om^\hdot  A )^\ka\br:=
((\Om^\hdot  A)\br)^\ka\sset 
(\Om^\hdot  A )\br$.

We recall 
that the cohomology group $H^n(\Om^{\,}\a,b)$ that occurs
in the  above displayed short exact sequence 
is isomorphic, as has been mentioned
in \S\ref{hoch_sect}, to the Hochschild homology
$H_n(A,A)$. Thus, Theorem \ref{hoch_thm} is an immediate
consequence of the  short exact sequence of the Proposition.

The following result, which was implicit
in  \cite{CQ2}, \cite{Lo}, will play an important role
in subsequent sections.

\begin{lem}\label{ka_iso} \vi The projection $(\Om^\hdot  A )\br\to\DR^\hdot  A $
restricts to a {\sf{bijection}}
$\dis (\Om^\hdot  A )^\ka\br\iso\DR^\hdot  A .$

\vii The map $b$ descends to 
a map $\dis b\br: (\Om^\hdot  A )\br\to
{\Om^{\idot-1} A}$.

\viii  The kernel of the map  $\dis b\br: (\Om^\hdot  A )^\ka\br\to
{\Om^{\idot-1} A},$
the restriction of $b\br$ to the space of
$\ka$-invariants, is isomorphic to $H^n(\Om^{\,}\a,b).$
\end{lem}

Both  Proposition \ref{bkappa} and Lemma \ref{ka_iso} will be proved
later,
in \S\ref{pf1}. 

\subsection{Special case: $H_1(A,A)$.}
 Note that the map $\ka$ acts trivially
on ${\Om^0 A}\cong A$. Therefore, on 1-forms,
 we have $\bi_\Delta=b$. Thus, using the
identification $H_1(A,A)=H^1(\omb, b)$,
 the short exact sequence of  Proposition \ref{bkappa}
reads
\beq{b1}
0\map H_1(A,A)\map \DR^1 A\,\stackrel{b=\bi_\Delta}\tooo\, [A,A]\map 0.
\eeq

The short exact  sequence \eqref{b1} may be obtained in an
alternate way as follows.
We apply  the right exact functor $(-)\br$
to the fundamental short exact sequence $0\to\Om^1 A\to
A\otimes A\to A\to0$, see \eqref{fund},
of $A$-bimodules.
The corresponding  long  exact sequence of Tor-groups
reads
$$\ldots\to H_1(A,A\otimes A)
\to H_1(A,A)\to(\Om^1 A)\br\to
(A\otimes A)\br\stackrel{c}\to A\br\to 0.
$$
Now, by definition of  Tor, one has
$H_k(A,A\otimes A)=0$ for all ${k>0}.$
Also, one has natural
 identifications $(\Om^1 A)\br=\DR^1 A,$
and
$(A\otimes A)\br\cong A.$
This way, the  map $c$ on the right of the displayed formula
may be identified with
the natural 
projection
$A\onto A/[A,A]$. Thus, $\Ker(c)=[A,A]$,
and the long  exact sequence above
reduces to the short exact sequence \eqref{b1}.

Further, it is immediate from definitions that
map $b=\bi_\Delta$ in \eqref{b1} is given
 by the 
formula 
$$
u\,dv\mapsto [u,v],\quad\text{hence}\quad
(\DR^1 A)\exact\sset
\Ker(\bi_\Delta)=H_1(A,A).
$$
This map $u\,dv\mapsto [u,v]$ has been considered
earlier by Cuntz and Quillen, see \cite{CQ1}.

\subsection{An application.} In this section we use
Theorem \ref{hoch_thm} in order to prove the following

\begin{prop}\label{ham_lemma}
Let  $A$ be an  algebra with  $H_2(A,A)=0$ and such that  the following
sequence
is exact (in which case $A$ may be called `{\em connected}')
\beq{exact}
0\map \k\map \DR^0 A\stackrel{d}\map\DR^1 A.
\eeq

Then, one has \smallskip

\npb{$H_1(A,A)=(\DR^1 A)\closed=(\DR^1 A)\exact$.}

\npb{There is a natural  vector space
isomorphism $(\DR^2 A)\closed\iso[A,A].$}
\end{prop}

\begin{proof} We will freely use the notation of 
\cite[\S4.1]{CBEG}. 

According to \cite{CBEG}, Proposition 4.1.4, for any algebra $A$
such that the sequence \eqref{exact} is exact, one has
the following  commutative diagram 
\beq{sqdiag}
\xymatrix{
\DR^1 A\ar[rr]^<>(0.5){d}\ar@{->>}[d]_<>(0.5){\bi_\Delta}&&
(\DR^2 A)\closed
\ar@{.>}[dll]_<>(0.5){\widetilde{\munc}}\ar[d]_<>(0.5){\bi_\Delta}\\
{[A,A]}^{}\ar@{^{(}->}[rr]^<>(0.5){d}&&[A,\Om^1 A].
}
\eeq
In this diagram, the map
  ${\widetilde{\munc}}$ is a lift of the
{\em noncommutative moment map}
 that has been introduced in \cite{CBEG}.

 Assuming that
$\HH_2(A)=0$, we deduce from
 the short exact sequence of
 Proposition \ref{bkappa} for $n=2$
 that the map
$\bi_\Delta: \DR^2A\to[A,\Om^1 A]$ is injective.

We now exploit diagram \eqref{sqdiag}.
The map  $d$ in the bottom row of the diagram is injective
due to acyclicity of the de Rham differential on $\omb$,
cf. eg. formula (2.5.1) in \cite{CBEG}. Further,
 the left vertical map $\bi_\Delta$ in the diagram
is surjective by \eqref{b1}.
Therefore, using commutativity of  diagram \eqref{sqdiag}
we deduce by diagram chase
that the upper horisontal map $d$ must be surjective,
and also the map  ${\widetilde{\munc}}$ must be bijective.

This yields both statements of Proposition
\ref{ham_lemma}.
\end{proof}

A version of
Proposition \ref{ham_lemma}  applies 
in the case where $A$ is the path algebra of a quiver
with $r$ vertices.
In that case, one has to consider algebras over a ground
ring $R:=\k\oplus\ldots\oplus\k$ ($r$ copies)
rather than over the base field $\k$.
The corresponding formalism has been worked out in \cite{CBEG}.

Write $\DR^\hdot_RA$ for the
relative 
 de Rham complex of $A$ and let $[A,\Om^1_RA]^R$ 
denote the vector space formed
by those elements of $[A,\Om^1_RA]$ which commute with $R$.

An analogue of Proposition \ref{ham_lemma}
implies the following result.

\begin{cor} Let $A$ be the path algebra of a quiver. Then
there is a natural  vector space
isomorphism $(\DR^2_RA)\closed\iso[A,A]^R.$\qed
\end{cor}

\pagebreak[3]
\section{Proofs}\label{pf_main}
\subsection{Proof of Lemma \ref{ka_iso} and Proposition \ref{bkappa}.}\label{pf1}
Our proof of Lemma \ref{ka_iso} follows
the proof of Lemma 2.6.8 in \cite{Lo}.
Write
$\Om^n:={\Om^n  A} $ and $\Om:=\oplus_n\,\Om^n$.

One observes  from definitions that 
$[A,\Om]=b\Om$ and $[dA,\Om]=(\Id-\ka)\Om.$
Hence,  we obtain, cf.~\cite{CQ1}:
$$[\Om,\Om]=[A,\Om]+[dA,\Om]=b\Om+(\Id-\ka)\Om.
$$
 We deduce 
that
$\Om\br=\Om/b\Om,$ and
$\dis\DR^\hdot A =\Om/[\Om,\Om]=\Om\br/(\Id-\ka)\Om\br.$
It follows in particular, since $b^2=0$, that $b$ descends to a well
defined map $b\br: \Om\br=\Om/b\Om\to\Om$.

Further, one has the following
standard identities, see \cite[\S2]{CQ2}
\beq{ident} \ka^n-\Id=b\ccirc\ka^n\ccirc d,
\quad\ka^{n+1}\ccirc d=d
\quad
\text{hold on}\en \Om^n,\quad \forall n=1,2,\ldots.
\eeq

The  Karoubi operator $\ka$ commutes with $b$,
hence induces a well-defined endomorphism of
 the vector space
 $\Om^n/b\Om^n,\, n=1,2,\ldots.$ Furthermore,
from the first identity in \eqref{ident}
we see that $\ka^n=\Id$ on
 $\Om^n/b\Om^n$.
Hence, we have a direct sum decomposition
$\Om\br=(\Om\br)^\ka\oplus(\Id-\ka)\Om\br$.
It follows that  the natural
projection $\Om\br=\Om/b\Om\onto \DR^\hdot A =\Om\br/(\Id-\ka)\Om\br$
restricts to an isomorphism
$(\Om\br)^\ka\iso\DR^\hdot A.$

Parts (ii) and (iii) of  Lemma  \ref{ka_iso} are
clear from the proof of Lemma 2.6.8 in \cite{Lo}. \qed

\begin{proof}[Proof of Proposition \ref{bkappa}.]
 The equality $\bi_\Delta=
b\ccirc (1+\ka+\ka^2+\ldots+\ka^{n-1}),$
in the statement of  the Proposition, is immediate from
the formula of Lemma~\ref{biDelta}(ii).

To prove the second statement we 
exploit the first identity in \eqref{ident}. Using
the formula for $\bi_\Delta$ and the fact
that $b$ commutes with $\ka$, we compute
$$(\ka-1)\ccirc\bi_\Delta=b\ccirc(\ka-1)\ccirc(1+\ka+\ka^2+\ldots+\ka^{n-1})
=b\ccirc(\ka^n-1)=b^2\ccirc\ka^n\ccirc d=0.
$$
Hence, we deduce that the image of $\bi_\Delta$ is contained
in $(b\Om)^\ka.$
Conversely, given any element $\al=b(\beta)\in (b\Om)^\ka,$
we find
$$\bi_\Delta(\beta)=
(1+\ka+\ka^2+\ldots+\ka^{n-1})\ccirc b(\beta)=n\cdot b(\beta)=n\cdot\al.
$$
Thus, we have
$\op{Im}(\bi_\Delta)=(b\Om)^\ka=(\,[A,\Om]\,)^\ka,$ since $b\Om={[A,\Om]}.$
Furthermore, it is clear that
the two maps
$(1+\ka+\ka^2+\ldots+\ka^{n-1})\ccirc b$
and $b$ coincide on
$(\Om^\hdot)^\ka\br$, hence,
have the same kernel.
The exact sequence of the Proposition  now
follows from Lemma \ref{ka_iso}.
\end{proof}

\subsection{Harmonic decomposition.} 
Our proof of Theorem \ref{main} is an adaptation of the strategy used in
\cite[\S2]{CQ2}, based on a  {\em Harmonic decomposition}
$$\bom=P\bom\oplus P^\perp\bom,\quad\text{where}\en
P\bom:=\Ker(\Id-\ka^2),\en P^\perp\bom:=\op{Im}(\Id-\ka^2).
$$

The differentials $B, b,$ and $d$ commute with $\ka$, hence
preserve Harmonic decomposition. Moreover, 
the differentials $B$ and $d$ are known to be
 proportional on $P\bom$. 
Specifically, exploiting  the second identity in
\eqref{ident}, it has been shown in
\cite{CQ2}, formula (11),
that on $P\bom^{_n}$
one has
\beq{agreeCQ}
B=(n+1)\cdot d\quad\text{on}\en P\bom^{_n},\quad\forall n=1,2,\ldots.
\eeq

Next,  we claim that
\beq{agree}
{\mathsf{(i)}}\en\bi_\Delta|_{P^\perp\bom}=0,
\quad\text{and}\quad
{\mathsf{(ii)}}\en\bi_\Delta=n\cdot b\quad\text{on}\en P\bom^{_n},\quad\forall n=1,2,\ldots.
\eeq

Claim (i) follows from Proposition \ref{bkappa}, since
$\op{Im}(\bi_\Delta)\sset \Ker(\Id-\ka)$
$\sset\Ker(\Id-\ka^2)=P\bom.$

To prove (ii), let $\al\in \bom^{_n}$. From the first  identity in
\eqref{ident}, we get $\al-\ka^n(\al)\in b\bom$.
Hence,  we find
$b\al-\ka^n(b\al)\in b^2\bom=0,$
 since $b^2=0$.
Thus, the operator $\ka$ has finite order on
$b\bom$, hence, also on $b(P\bom)$.
 But, for any operator $T$ of finite order,
one has $\Ker(\Id-T)=\Ker((\Id-T)^2).$
It follows  that if $\al\in P\bom^{_n}$,
then we have $b\al\in\Ker((\Id-\ka)^2)=\Ker(\Id-\ka).$
We conclude that the element $b\al$ is
fixed by $\ka$. Therefore, by Proposition \ref{bkappa},
we obtain
$$\bi_\Delta(\al)=(1+\ka+\ka^2+\ldots+\ka^{n-1})\ccirc b(\al)=
n\cdot b(\al),$$
and \eqref{agree} is proved. 

\subsection{Proof of Theorem \ref{main}.}\label{pf2}
Harmonic decomposition being stable under
all four differentials $B, b, d,$ and $\bi_\Delta$,
we may analyze the homology of each of the direct summands,
$P\bom$ and $P^\perp\bom$, separately. 

First of all,  it has been shown by Cuntz-Quillen, see
Proposition 4.1(1) in  \cite{CQ2}, that the complex
 $(P^\perp\bom[t], B+t\cdot b)$
is  acyclic.

Further,  the complex
$(\bom,d)$  being acyclic, cf. \cite[\S1]{CQ2} or \cite{CBEG} formula (2.5.1),
we deduce
\beq{acy}
\text{\em Each of the complexes}\en
(P\bom,d)\en \text{\em and}  \en (P^\perp\bom,d)\en
\text{\em is acyclic.}
\eeq

Now,
the map $\bi_\Delta$ vanishes on
$P^\perp\bom$
by \eqref{agree}(i). Hence, on $P^\perp\bom[t]$ we have
$d+t\cdot \bi_\Delta=d.$
 Therefore, we conclude using \eqref{acy}
that  $(P^\perp\bom[t], d)$, hence, $(P^\perp\bom[t], d+t\cdot
 \bi_\Delta)$
is an acyclic complex.

Thus, to complete the proof of the theorem, 
it suffices to construct an isomorphism
of complexes $\dis F^\hdot: (P\bom[t], d+t\cdot \bi_\Delta\bigr)\iso
(P\bom[t],  B+t\cdot b)$. Such an isomorphism is provided
by the assignment $F^{(n)}: \al\mapsto n!\cdot\al,$
$\forall\al\in(P\bom)^n,$
since for any $n=0,1,\ldots$
we clearly have 
$$B\ccirc F^{(n)}=F^{(n+1)}\ccirc d,\quad\text{and}\quad
b\ccirc F^{(n)}=F^{(n-1)}\ccirc \bi_\Delta.
\qquad\Box
$$

\section{Miscellaneous results}
\subsection{First order deformations based of free products.}\label{defor1}
Below, we outline a new version of deformation theory of
associative
algebras, based on the free product construction
$A*\k[t]$.
We are going to study  {\em first order} deformations
of an associative
algebra $A,$
 in which the  deformation
parameter $t$ is {\em not} assumed to be central.

Motivated by the classical approach,
we first identify the vector space
$A_t/(A_t^+)^2$ with
$A\oplus (A\o A),$ using \eqref{AtA}.
Thus,   we are interested in associative products
on the vector space $A\oplus (A\o A)$ that have
the following form 
\beq{first}
\bigl(u\oplus(u'\o u'')\bigr)\; \times\; 
\bigl(v\oplus(v'\o v'')\bigr)
\stackrel{\starb}\mto
uv\oplus\bigl(u'\o u''v +uv'\o v'' +\be(u,v)\bigr),
\eeq
where $\be: A\times A \to A\o A$ is a certain
$\k$-bilinear map.
 
These products are taken up to an equivalence.
Specifically, for any $\k$-linear map  $f: A\to A\o A$,
one defines a linear bijection
$$\tilde{f}:\
A\oplus (A\o A)\too A\oplus (A\o A),\quad
 u\oplus(u'\o u'')\mto u\oplus(u'\o u''+ f(u)).$$
Given a product $\starb$ and
a map $f$, one defines a new
 product by transporting the structure via $\tilde{f}$,
that is, by the formula
$\dis
x\star_{_\gamma}y:=\tilde{f}\inv(\tilde{f}(x)\starb\tilde{f}(y)).$
We say that the products $\star_{_\gamma}$ and $\starb$ are
 {\em equivalent}.

 Similarly to  conventional theory, 
a   classification of equivalence classes of associative products
may be given
in terms of Hochschild cohomology. 
Specifically, one performs
a routine calculation, completely analogous
to the  classical one  due to Gerstenhaber, that yields the following

\npb{A first product $\starb$ as in \eqref{first}
is associative $\Longleftrightarrow$
$\be\in C^2(A, A\o A)$ is a Hochschild $2$-cocycle with coefficients
in $A\o A$.}

\npb{Twe products $\starb$ and $\star_{_\gamma}$ corresponding to  2-cochains $\be$ and
$\gamma$ are equivalent $\Longleftrightarrow$
$\be-\gamma$ is a  Hochschild coboundary.}
\medskip

\noindent
Thus, we deduce

\begin{prop} Equivalence classes of 
 associative  products, as  in \eqref{first},  are in one-to-one correspondence
with the elements of $H^2(A,A\o A)$, the second  Hochschild cohomology
group
of the $A$-bimodule $A\o A$.\qed
\end{prop}

\subsection{Another construction of the extended de Rham complex.}
\label{another}
Let $A$ be an algebra and $M$ an $A$-bimodule.
We consider
 ${A^{\o 2}}\oplus M$ as an $A$-bimodule and introduce
an $A$-bimodule map
$\pi: {A^{\o 2}}\oplus M\to A,\, (a'\o a'')\oplus u
\mapsto a'a''.$

Further, we equip  the bimodule ${A^{\o 2}}\oplus M$ with
a grading in two different ways. The first grading
assigns ${A^{\o 2}}$ grade degree 2, while
the second  grading
assigns ${A^{\o 2}}$ grade degree $-1$. In both cases,
the direct summand $M$ is placed in degree $0$.
These two gradings on  ${A^{\o 2}}\oplus M$ induce the corresponding
gradings on the tensor algebra
$T_A({A^{\o 2}}\oplus M),$ to be denoted
$T^t_A({A^{\o 2}}\oplus M)$ and
$T^\tau_A({A^{\o 2}}\oplus M),$ respectively

The following result is a generalization
of Proposition \ref{drin}.
\begin{prop}\label{tsy} \vi For any $A$-bimodule $M$, the map
$\dis
(a'\o a'')\oplus m\mto a'ta''+m
$
can be uniquely extended, by multiplicativity, to a graded algebra isomorphism
$\dis T^t_A({A^{\o 2}}\oplus M)\iso (T_AM)_t$, resp.,
$\dis T^\tau_A({A^{\o 2}}\oplus M)\iso (T_AM)_\tau.$

\vii The  differential 
$\frac{d}{d\tau}:   (T_AM)_\tau\to (T_AM)_\tau,$ cf. \S\ref{ttau}, is transported
under the isomorphism in \vi to a  degree 1 map
$\pa : 
T^\tau_A({A^{\o 2}}\oplus M)\to
T^\tau_A({A^{\o 2}}\oplus M)$ given,
 for any homogeneous $x_1,\ldots,x_n\in {A^{\o 2}}\oplus M,$
by 
\begin{align}\label{tsy2}
\pa (x_1\o\ldots\o x_n)=&\sum_{k=1}^{n-1} 
(-1)^{|x_1,\ldots,x_{k-1}|}
x_1\o\ldots\o x_{k-1}\o\pi(x_k)x_{k+1}\o\ldots\o x_n,\nonumber\\
&|x_1,\ldots,x_{k-1}|:=k-1+\deg x_1+\ldots+\deg x_{k-1}.
\end{align}
\end{prop}

\begin{proof} The assignment 
$(a'\o a'')\oplus m\mapsto a'ta''+m,$ clearly gives an
$A$-bimodule map $A^{\o 2}\oplus M\to
(T_AM)_t$. This map can be extended, by the universal
properties,
to an algebra homomorphism $T_A({A^{\o 2}}\oplus M)\to
(T_AM)_t$. To show that this  homomorphism is an isomorphism,
we construct an inverse map as follows.

The algebra  $(T_AM)_t=(T_AM)*\k[t]$ is freely generated by 
the element $t$ and the subalgebra $T_AM\sset(T_AM)_t.$
Therefore, constructing an algebra homomorphism
$(T_AM)_t\to T_A({A^{\o 2}}\oplus M)$ amounts to
giving an algebra homomorphism
$f: T_AM\to T_A({A^{\o 2}}\oplus M)$ and an additional
 element
$f(t)\in T_A({A^{\o 2}}\oplus M)$,
 chosen arbitrarily.
Thus, we may define the homomorphism $\dis f: T_AM\to T_A({A^{\o 2}}\oplus M)$
to be the natural algebra imbedding
$\dis T_AM\into T_A({A^{\o 2}}\oplus M),$
induced by the $A$-bimodule imbedding
$\dis M=0\oplus M\into{A^{\o 2}}\oplus M.$
Further,  we set $\dis f(t):=1_A\o 1_A\in A^{\o 2}\sset
T^1_A({A^{\o 2}}\oplus M)$.
It is straightforward to check that the resulting
homomorphism
$\dis (T_AM)_t\to T_A({A^{\o 2}}\oplus M)$ is indeed an inverse
of the homomorphism in the opposite direction constructed
in the first paragraph of the proof. 

This completes the
proof of part (i) (compatibility with the corresponding
 gradings is clear). Part (ii) is proved
by a straightforward computation.
\end{proof}

In the special case $M=0$, part (i) of the above Proposition
 is due to Van den Bergh
(private communication). 
The differential $\pa $ in \eqref{tsy2} is a special
case of the generalized bar differential introduced by 
B. Tsygan \cite{T}.

Applying Proposition \ref{tsy} in the special
case $M=\oma,$ we obtain 
\begin{cor}\label{big} There is a canonical 
algebra isomorphism
$$
\Om_t A\,\bigl(=(\Om^{\,}\a)_t\bigr)\cong T^t_A(A^{\o 2}\oplus\oma).
$$

For any double derivation $\Th\in\dder\a$, the 
degree $-1$ graded derivation
 $T^t_A(A^{\o 2}\oplus\oma)\to T^t_A(A^{\o 2}\oplus\oma),$
corresponding to the contraction
operation $i_\Th:\Om_t A\to\Om_t A,$
is induced by the  following $A$-bimodule map
$$A^{\o 2}\oplus\oma\map A^{\o 2}\oplus\oma,\quad
(a'\o a'')\oplus \al
\mto (i'_\Th\al\o i''_\Th\al)\oplus 0.
$$
\end{cor}

{\small{

}}

\smallskip

{\small
\noindent
{{Department of Mathematics, University of Chicago,
Chicago, IL
60637, USA;}\newline
E-mail address:\; {\tt ginzburg@math.uchicago.edu}}

\end{document}

We give a new construction of cyclic homology of an associative 
algebra A that does not involve Connes' differential. Our approach
is based on an extended version of the complex \Omega A, of noncommutative
differential forms on A, and is similar in spirit to the de Rham 
approach to equivariant cohomology.  Indeed, our extended complex
maps naturally to the equivariant de Rham complex of any representation scheme Rep A.

We define cyclic homology as the cohomology of the total complex 
(\Omega A)[t], d+t \cdot i, arising from two anti-commuting
differentials, d and i, on \Omega A. The  differential d, that replaces 
the Connes differential B, is the Karoubi-de Rham differential. The  
differential i that replaces the Hochschild differential b, is a  map  analogous to 
contraction with a vector field. This new map  has no commutative counterpart.

 User-ID: math.KT/0505236
 Password: 8nf2z